\documentclass[]{article}     
\usepackage{amsthm}
\usepackage{authblk}
\usepackage[margin=2cm]{geometry}

\RequirePackage{amsmath}

\RequirePackage{amssymb}
\RequirePackage{amsfonts}

\RequirePackage{bbm}

\RequirePackage{graphicx}

\RequirePackage{booktabs}

\RequirePackage{url}

\RequirePackage{hyperref}

\RequirePackage{xcolor}
\RequirePackage{color}

\RequirePackage{wasysym}
\RequirePackage{listings}
\RequirePackage{caption}
\RequirePackage{subcaption}

\RequirePackage{stmaryrd}

\RequirePackage{lineno}

\RequirePackage{xparse}

\RequirePackage{float}

\RequirePackage{varioref}

\RequirePackage{multirow}

\RequirePackage{cellspace}

\RequirePackage{shellesc}

\RequirePackage[ruled,vlined,linesnumbered]{algorithm2e}

\RequirePackage{tikz}
\usetikzlibrary{automata, positioning, arrows, external, matrix}

\RequirePackage{pgfplots}
\usepgfplotslibrary{fillbetween}
\usepgfplotslibrary{groupplots}
\usepgfplotslibrary{external}

\DeclareMathAlphabet{\pazocal}{OMS}{zplm}{m}{n}

\providecommand{\bR}{\mathbb{R}}

\providecommand{\pV}{\pazocal{V}}
\providecommand{\pX}{\pazocal{X}}

\providecommand{\pT}{\pazocal{T}}

\providecommand{\pM}{\pazocal{M}}
\providecommand{\pN}{\pazocal{N}}
\providecommand{\pP}{\pazocal{P}}

\providecommand{\bs}{\boldsymbol}
\providecommand{\bx}{\bs{x}}
\providecommand{\bu}{\bs{u}}
\providecommand{\bv}{\bs{v}}
\providecommand{\bn}{\bs{n}}
\providecommand{\bA}{\bs{A}}

\providecommand{\half}{\frac{1}{2}}

\providecommand{\vertiii}[1]{{\left\vert\kern-0.15ex\left\vert\kern-0.15ex\left\vert #1
    \right\vert\kern-0.15ex\right\vert\kern-0.15ex\right\vert}}

\NewDocumentCommand{\curlii}{sO{}m}
{
	\IfBooleanTF{#1}
    {\dgalext{#3}}
    {\dgalx[#2]{#3}}
}

\NewDocumentCommand{\dgalext}{m}{  \sbox0{    \mathsurround=0pt     $\left\{\vphantom{#1}\right.\kern-\nulldelimiterspace$  }  \sbox2{\{}  \ifdim\ht0=\ht2
    \{\kern-.625\wd2 \{#1\}\kern-.625\wd2 \}  \else
    \left\{\kern-.7\wd0\left\{#1\right\}\kern-.7\wd0\right\}  \fi
}

\NewDocumentCommand{\dgalx}{om}{  \sbox0{\mathsurround=0pt$#1\{$}  \sbox2{\{}  \ifdim\ht0=\ht2
    \{\kern-.625\wd2 \{#2\}\kern-.625\wd2 \}  \else
    \mathopen{#1\{\kern-.7\wd0 #1\{}
    #2
    \mathclose{#1\}\kern-.7\wd0 #1\}}
  \fi
}

\tikzset{
  partial ellipse/.style args={#1:#2:#3}{
    insert path={+ (#1:#3) arc (#1:#2:#3)}
  }
}

 \graphicspath{{./pngfigs/}}
\providecommand{\keywords}[1]
{
  \small	
  \textbf{\textit{\ \ \ \ Keywords:}} #1
}

\makeatletter
\if@twocolumn
\newcommand{\whencolumns}[2]{
#2
}
\else
\newcommand{\whencolumns}[2]{
#1
}
\fi
\makeatother

\newtheorem{remark}{Remark}

\begin{document}

\title{A Multigrid Method for a Nitsche-based Extended Finite Element Method}

\author[]{Hardik Kothari\thanks{hardik.kothari@usi.ch}}
\author[]{Rolf Krause\thanks{rolf.krause@usi.ch}}
\affil[]{Institute of Computational Science, Universit\`{a} della Svizzera italiana, Lugano, Switzerland}

\maketitle

\begin{abstract}
   We present a tailored multigrid method for linear problems stemming from a Nitsche-based extended finite element method (XFEM).
Our multigrid method is robust with respect to highly varying coefficients and the number of interfaces in a domain.
It shows level independent convergence rates when applied to different variants of Nitsche's method.
Generally, multigrid methods require a hierarchy of finite element (FE) spaces which can be created geometrically using a hierarchy of nested meshes.
However, in the XFEM framework, standard multigrid methods might demonstrate poor convergence properties if the hierarchy of FE spaces employed is not nested.
We design a prolongation operator for the multigrid method in such a way that it can accommodate the discontinuities across the interfaces in the XFEM framework and recursively induces a nested FE space hierarchy.
The prolongation operator is constructed using so-called pseudo-$L^2$-projections; as common, the adjoint of the prolongation operator is employed as the restriction operator.

The stabilization parameter in Nitsche's method plays an important role in imposing interface conditions and also affects the condition number of the linear systems. We discuss the requirements on the stabilization parameter to ensure coercivity and review selected strategies from the literature which are used to implicitly estimate the stabilization parameter.
Eventually, we compare the impact of different variations of Nitsche's method on discretization errors and condition number of the linear systems.
We demonstrate the robustness of our multigrid method with respect to varying coefficients and the number of interfaces and compare it with other preconditioners.
 \end{abstract}

\section{Introduction}
In the last two decades, unfitted finite element methods have become quite popular.
An unfitted method can be defined as any method where the computational domain does not match the mesh exactly.
These methods ease the strict requirement of traditional finite element methods (FEM), where a mesh has to be generated to represent the computational domain.
In many cases, it can be a demanding task to create high-quality meshes for complex geometries and failing to do so can result in usually sub-optimal approximation properties of FEM.
Geometrically unfitted methods just require a background mesh, and a finite element space defined on the background mesh.
Clearly, the latter has to be modified to enforce the boundary conditions or the interface conditions.
Here, an interface can be described as codimension one entity embedded in the domain, across which a function exhibits non-smooth properties.

There is a huge variety of unfitted methods. The fictitious domain method can be listed as one of the oldest variants of an unfitted method~\cite{glowinski_fictitious_1994}.
The eXtended Finite Element Method (XFEM) has been introduced for problems in fracture mechanics with crack propagation~\cite{fries_extended/generalized_2010,sukumar_extended_2000,dolbow_extended_2001}.
The XFEM method reduces the computational burden of remeshing but retains the robustness of the FEM~\cite{belytschko_elastic_1999}.
A similar XFEM approach was taken for the two-phase flow problems in fluid dynamics to enrich the pressure variable~\cite{gros_finite_2006,reusken_analysis_2008,SvenGross2011-04-25}.
In both methods, the crack and the interface between the two-phases evolve over time.
In the traditional fitted FEM repeated remeshing of the domain is required to accomodate the evolving interfaces, while in the XFEM framework this problem can be dealt with by enriching the standard finite element (FE) spaces.
The Finite Cell Method (FCM) can be considered as an extension of the fictitious domain method on higher-order function space~\cite{parvizian_finite_2007,schillinger_finite_2015}.
The FCM has been recently extended by employing the spline-based finite elements for harmonic and bi-harmonic problems~\cite{embar_imposing_2010}.
Other examples of unfitted methods can be given as the CutFEM method~\cite{burman_cutfem:_2015,claus_cutfem_2015}, the immersed boundary methods~\cite{peskin_immersed_2002}, the trace finite element method~\cite{olshanskii_trace_2017}, etc.

In all variants of unfitted methods, the boundaries or the interfaces are not resolved by the mesh explicitly.
In such cases, it is not possible to enforce boundary conditions (BC) or interface conditions explicitly as nodal values.
Dirichlet BC/interface conditions are enforced by using either the penalty method~\cite{babuska_finite_1973}, the method of Lagrange multipliers~\cite{babuska_finite_1973-1} or Nitsche's method~\cite{nitsche_uber_1971}.
Although the penalty method is trivial to implement, it is not widely used.
The penalty method is variationally inconsistent and it can produce a highly ill-conditioned system if a large penalty parameter is required to ensure optimal convergence of discretization error.
The method of Lagrange multipliers is an attractive option, but the method leads to optimal convergence of discretization error only if a proper finite element space is chosen for the boundaries and/or interfaces~\cite{ji_strategies_2004}.
The choice of a finite element space for the Lagrange multiplier is not obvious as the method is not necessarily stable if the discrete inf-sup condition is not satisfied~\cite{bechet_stable_2009}.
The method of Lagrange multipliers requires a mixed-finite element formulation and the system of linear equations has a saddle-point structure, which adds additional computational complexity for the solution methods.
Nitsche's method is an alternative to these methods to enforce Dirichlet BC or interface conditions.
The method can be regarded as a variationally consistent penalty method.
But a different interpretation of the method is as a stabilized Lagrange multiplier method, where the Lagrange multiplier is explicitly expressed by its physical interpretation in the primal variable~\cite{stenberg_techniques_1995}.
The method is utilized in many different discretization methods because of its versatility.
In the discontinuous Galerkin (DG) method, Nitsche's method is used to enforce the continuity between each element faces~\cite{arnold_interior_1982}.
The method is also used in domain decomposition method to mortar the interfaces between non-matching meshes~\cite{brezzi_stabilization_1997,becker_finite_2003}.
Nitsche's method is also a popular choice to enforce boundary conditions for mesh-free methods and particle methods~\cite{griebel_particle-partition_2003}.

In the context of XFEM, Nitsche's method is introduced and analyzed for elliptic interface problems with the discontinuous coefficients for unfitted meshes~\cite{hansbo_unfitted_2002}.
The extension of the method is also provided for strong and weak discontinuities in solid mechanics and the optimal convergence properties of Nitsche's method in the XFEM framework is shown~\cite{hansbo_finite_2004}.
After the initial works~\cite{hansbo_unfitted_2002,hansbo_finite_2004}, many efforts have been made to improve the robustness of Nitsche's formulation in the XFEM framework.
In this work, we discuss selected strategies for estimating the stabilization parameter~\cite{ruess_weakly_2013,lehrenfeld_removing_2016} and discuss the ghost penalty stabilization method~\cite{burman_ghost_2010}.

In the unfitted methods, a background mesh captures the computational domain of arbitrary shape, thus the elements are allowed to cut arbitrarily by the boundaries or interfaces.
This could give rise to a highly ill-conditioned system of linear equations.
Due to this reason, it becomes essential to develop efficient solution strategies for such discretization methods.
In this work, we propose a tailored multigrid method for solving a system of equations arising from Nitsche-based XFEM.

Multilevel methods are ideal iterative solvers for many large-scale linear/nonlinear problems as they are of optimal complexity~\cite{WolfgangHackbusch1986-01-14}.
The optimal complexity implies that the convergence rate of the multilevel methods is bounded from above independently from the size of the problem.
The robustness of multilevel iteration results from a sophisticated combination of smoothing iterations and coarse level corrections.
Ideally, these components are complementary to each other as they reduce error in different parts of the spectrum.
Traditionally, the mesh hierarchy for multilevel methods is created by either coarsening or refinement strategies, and simple interpolation operator and its adjoint are used to transfer the information between different levels.

There have been some efforts to develop multilevel solution strategies for the XFEM discretization.
Initial approaches propose to modify the algebraic multigrid method (AMG).
A domain decomposition-based AMG preconditioner is proposed for the fracture problems~\cite{berger-vergiat_inexact_2012}, where the domain is decomposed into `cracked' and `intact' domain and AMG is applied to the `intact' domain, and the `cracked' domain is solved with the direct solvers.
In an alternative approach, known as a quasi-algebraic multigrid method, the sparsity pattern of the interpolation operator is modified to prevent the interpolation across the interfaces~\cite{hiriyur_quasi-algebraic_2012}.
Recently, a new multigrid method is also proposed for the elliptic interface problems, with an interface smoother~\cite{ludescher_multigrid_2018}.

In this work, we propose a multigrid method specifically tailored for the unfitted finite element methods.
This method borrows the ideas from the semi-geometric multigrid (SMG) method.
The SMG method has been developed to overcome the shortcomings of the standard multilevel methods~\cite{dickopf_multilevel_2010,dickopf_efficient_2009}.
This method facilitates the creation of transfer operators between non-nested mesh hierarchy.
Hence the non-nested meshes can be created independently for complex geometry and the transfer operators are computed using $L^2$-projections.
Conventionally, unfitted methods use structured background meshes, but still, due to the arbitrary location of the interfaces and/or the boundaries, a non-nested mesh hierarchy emerges.
We exploit the same strategy to compute the transfer operator by means of $L^2$-projection between the successive meshes in the hierarchy for unfitted methods.

The outline of this paper is given as follows.
We introduce a model problem and the framework for the XFEM discretization in section~\ref{sect:model}.
We discuss Nitsche's method and different approaches to estimate the stabilization parameter.
In section~\ref{sect:multigrid}, we introduce the framework of our tailored multigrid method for XFEM.
We present a strategy to create a hierarchy of nested FE spaces from the hierarchy of non-nested XFEM meshes by means of $L^2$-projections. 
We discuss the $L^2$-projection and pseudo-$L^2$-projections for computation of the transfer operators and extend the same idea in the XFEM framework.
Lastly, in section~\ref{sect:results} we show the results of the numerical experiments.
We compare the effect of the different variations of Nitsche's methods on the condition number of the linear system and the discretization errors.
We also compare the SMG method with other preconditioners and demonstrate the robustness of our multigrid method with respect to highly varying coefficients and number of interfaces.
 \section{A Model Problem} \label{sect:model}
In this section, we consider the Poisson problem with discontinuous coefficients and provide a framework for XFEM discretizations.
We discuss Nitsche's method for imposing interface conditions and explore a few strategies to make the method stable.

We consider a Lipschitz domain $\Omega \subset \bR^d$, $d = 2,3,$ with interface $\Gamma$  which decomposes the domain, $\Omega$, into two non-overlapping subdomains $\Omega_1$ and $\Omega_2$, such that $\Omega = \Omega_1 \cup \Omega_2 \cup \Gamma$.
The interface is defined as $\overline{\Gamma} = \partial\Omega_1 \cap \partial\Omega_2$ and is assumed to be sufficiently smooth.
For simplicity, the interface $\Gamma$ is defined as polygonal. We define a sufficiently regular function $u_i:\Omega_i \cup \Gamma \to\bR$ as a pair $\{u_1,u_2\}=:u$ in $\Omega$.
The jump of $u$ on the interface is defined as
\[
  \llbracket u\rrbracket := u_1\vert_\Gamma-u_2\vert_\Gamma,
  \label{eq:jump}
\]
where $u_i\vert_{\Gamma}$ is the restriction of $u_i$ to $\Gamma$.

We consider a stationary diffusion problem with discontinuous coefficient $\alpha$ as
\begin{equation}
  \begin{aligned}
    -\nabla \cdot \alpha \nabla u             & = f &  & \text{in } \Omega_1 \cup \Omega_2, \\
    u                                         & = 0 &  & \text{on } \partial \Omega,        \\
    \llbracket u \rrbracket                   & = 0 &  & \text{on } \Gamma,                 \\
    \llbracket \alpha\nabla_{\bn} u\rrbracket & = 0 &  & \text{on } \Gamma,
  \end{aligned}
  \label{eq:dd_poisson_eq}
\end{equation}
where $f \in L^2(\Omega)$.
The coefficient $\alpha \in \bR^+$ is piecewise constant defined as
\[
  \alpha(\bx) = \alpha_i \geqslant \alpha_0 >0\quad \forall \bx \in \Omega_i.
\]
Here, $L^2$ defines a Lebesgue space of square integrable function on the domain $\Omega$, with inner product $(a,b)_{L^2(\Omega)}: = \int_\Omega a b \ d\Omega$, and the induced $L^2$-norm, $\|\cdot\|^2_{L^2(\Omega)} = {(\cdot,\cdot)_{L^2(\Omega)}} $.
Additionally, $H^k$ is Sobolev space with function and its $k^{th}$ weak derivative in $L^2$-space, and the associated norm is denoted by $\|\cdot\|_{H^k(\Omega)}$.

In problem~\eqref{eq:dd_poisson_eq}, continuity of the function $u$ and the continuity of the flux across the interface is enforced.
Also, we define the outward flux from the interface as $\nabla_{\bn} u = \bn \cdot \nabla u$.
For definiteness, we take the unit normal $\bs{n}$ as the outward pointing normal to $\Omega_1$ on $\Gamma$, $\bs{n}=\bs{n}_1=-\bs{n}_2$.

Problem~\eqref{eq:dd_poisson_eq} is consistent with the Poisson problem for the continuous coefficients, $\alpha_1 = \alpha_2$.
For the standard FEM approach, the interface conditions in~\eqref{eq:dd_poisson_eq} can only be imposed if nodes are placed on the interface explicitly.
For the fitted finite element method, the detailed analysis of the interface problem with discontinuous coefficients is carried out under the assumption that the function is only non-smooth at the interface, the problem has a unique solution in $H^2$ in each convex subdomain $\Omega_i$~\cite{chen_finite_1998}.

\subsection{XFEM Discretization}

In standard FEM discretization, the interface across which the function $u$ is discontinuous has to be aligned with the element faces.
In contrast, for the XFEM approach this requirement is relaxed, and the interface is allowed to be anywhere in the domain.
The XFEM discretization captures the interfaces by enriching the FEM solution space by duplicating the elements which are intersected by the interface.
The new degrees of freedom (dofs) are then associated with the duplicated elements.

\begin{figure*}[]
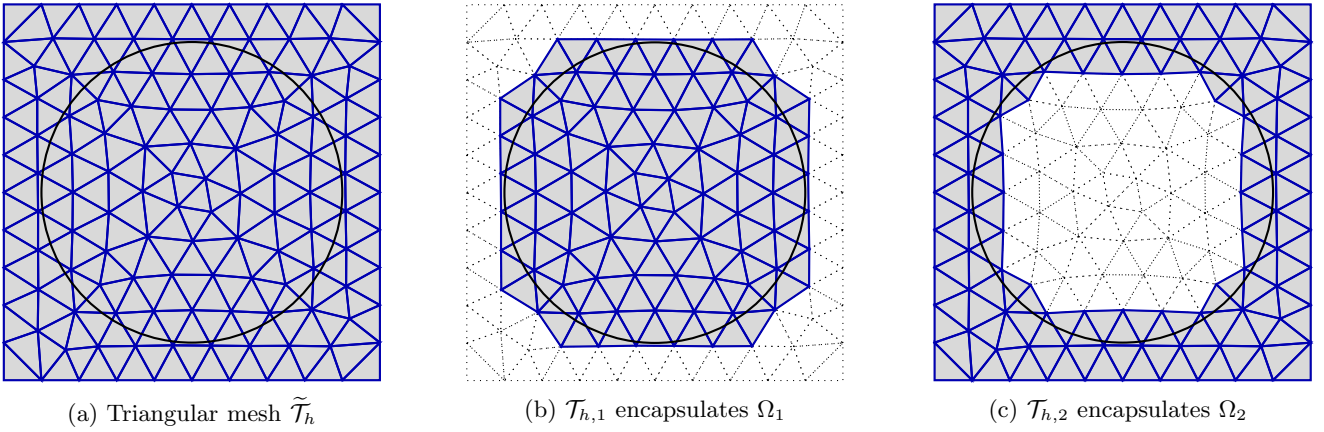

  \centering
  \subcaptionbox{Triangular mesh $\widetilde{\pT}_h$}[.3\textwidth]{

  }
  \hfill
  \subcaptionbox{$\pT_{h,1}$ encapsulates $\Omega_1$}[.3\textwidth]{%
  }
  \hfill
  \subcaptionbox{$\pT_{h,2}$ encapsulates $\Omega_2$}[.3\textwidth]{
    %
  }
  \caption[2D mesh]{2D Triangular mesh (a) is divided into two parts $\pT_{h,1}$ and $\pT_{h,2}$ as depicted in (b) and (c) due to the interface $\Gamma$, elements which belong to the respective domain are shaded}
  \label{fig:domain_decomp}
\end{figure*}

We assume a shape regular, quasi-uniform, conforming triangulation $\widetilde{\pT}_h$ on the domain $\Omega$.
We use $\widetilde{\pT}_h$ as background triangulation which captures both subdomains $(\Omega_1 \cup \Omega_2)\subseteq \widetilde{\pT}_h$.
Let $h_K$ be the diameter of the element $K$, and $h=\max_{K\in \widetilde{\pT}_h} h_K$.
We define the mesh associated with each subdomain as
\[
  \pT_{h,i} = \{ K\in \widetilde{\pT}_h: K\cap \Omega_i \neq \emptyset\},
\]
each subdomain $\Omega_i$ is encapsulated by the background mesh, $\Omega_i \subset \pT_{h,i}$, as shown in Figure \ref{fig:domain_decomp}.
Additionally, we define a set of elements that are intersected with the interface $\Gamma$ by
\[
  \pT_{h,\Gamma}=\{K \in \widetilde{\pT}_h: K\cap \Gamma \neq \emptyset \}.
\]
The interface triangulation $\pT_{h,\Gamma}$ is doubled: it is denoted as $\pT^i_{h,\Gamma} \subset \pT_{h,i}$ and it is part of both encapsulating meshes.
Additionally, for any element $K$, let $K_i=K\cap\Omega_i$ be part of $K$ in domain $\Omega_i$ and for $K\in \pT^i_{h,\Gamma}$, let $\Gamma_K:=\Gamma \cap K$ be part of $\Gamma$ in $K$.

We define continuous low order finite element space over the triangulation $\widetilde{\pT}_h$, which vanishes on the boundary as
\[
  \widetilde{\pV}_h = \{v\in H^1(\Omega) : v\vert_{K} \in \pP_1,\ v\vert_{\partial \Omega}=0,\ \forall K\in \widetilde{\pT}_h \},
\]
where $\pP_1$ denotes the space of piecewise linear functions.
Following the original XFEM works~\cite{belytschko_elastic_1999,moes_finite_1999}, we use the Heaviside function
\begin{equation}
  \label{eq:Heaviside}
  \pazocal{H}_i:\bR^d \to \bR, \qquad \pazocal{H}_i(\bx) = \begin{cases}
    1, & \quad \forall \bx \in \Omega_i,                  \\
    0, & \quad \forall \bx \in \Omega \setminus \Omega_i.
  \end{cases}
\end{equation}
The Heaviside function is used to restrict the support of the finite element space to domain $\Omega_i$, thus
the space of finite elements in the domain $\Omega_i$ is defined as
\[
  \pV_{h,i} = \pazocal{H}_i(\bx)\widetilde{\pV}_h.
\]
We seek the approximation $u_h := u_{h,1} \oplus u_{h,2}$ in space $\pV_h := \pV_{h,1} \oplus \pV_{h,2}$.
From the definition of the FE space, it is clear that the function is allowed to be discontinuous across the interface.
For more clarity, we show how the basis are defined on each space.
Let $\widetilde{\Phi}_h = (\widetilde{\phi}^p_h)_{p\in \widetilde{\mathcal{N}}_h}$ be nodal basis of $\widetilde{\pV}_h$ where $\widetilde{\mathcal{N}}_h$ denotes the set of nodes of the background triangulation.
Then, we obtain the ``cut" basis $\Phi_{h,i} = ({\phi}^p_h)_{p\in\mathcal{N}_{h,i}}, i=1,2$, with
\[
  \mathcal{N}_{h,i} = \{p\in \widetilde{\mathcal{N}}_h : \text{supp}(\phi^p_h)\cap \Omega_i \neq \emptyset \},
\]
and
\[
  \phi^p_h = \mathcal{H}_i(\bx) \widetilde{\phi}^p_h, \qquad \forall p \in \mathcal{N}_{h,i},
\]
where $\pazocal{N}_{h,i}$, denotes the set of nodes of the mesh associated with subdomain $\Omega_i$.
Following the same strategy, we define the span of the nodal basis function of the FE space $\pV_h$ as, $\Phi_h := \Phi_{h,1}\oplus \Phi_{h,2}$ and the set of nodes associated with the mesh $\pT_h$ is given by $\pN_h := \pN_{h,1} \oplus \pN_{h,2}$.
\subsection{Nitsche's Formulation}
In this section, we discuss Nitsche's formulation for the XFEM discretization denoted by~\eqref{eq:NXFEM}~\cite{hansbo_unfitted_2002}.
In this formulation, the interface condition is enforced weakly with a variant of Nitsche's method.
Before moving to the weak formulation, we define the average function as
\begin{equation*}
  \curlii[\big]{u} = \big( \beta_1u_1 + \beta_2u_2\big)
\end{equation*}
and
\begin{equation*}
  \curlii[\big]{\alpha \nabla_{\bn} u } = \big(\beta_1 \alpha_1 \nabla_{\bn} u_1 + \beta_2 \alpha_2 \nabla_{\bn} u_2 \big),
\end{equation*}
where $\beta_i \in \bR^+$ are the weighting parameters.
The abstract variational formulation of the problem~\eqref{eq:dd_poisson_eq} in the context of XFEM can be written as follows: find $u_h \in\pV_h$, such that
\begin{equation}
  \label{eq:NXFEM}
  a(u_h,v_h) = F(v_h), \qquad \forall v_h \in \pV_h,\tag{N}
\end{equation}
where
\begin{equation}
  \whencolumns
  {
    a(u_h,v_h)=\sum_{i=1}^2\Bigg(
    \int\limits_{\Omega_i} \alpha \nabla u_h \cdot \nabla v_h \ d\Omega_i \Bigg)
    - \int\limits_{\Gamma} \curlii[]{\alpha \nabla_{\bn} u_h} \llbracket v_h \rrbracket \ d\Gamma
    - \int\limits_{\Gamma} \llbracket u_h \rrbracket \curlii[]{\alpha \nabla_{\bn} v_h} \ d\Gamma
    + \int\limits_{\Gamma} \gamma \llbracket u_h \rrbracket \llbracket v_h \rrbracket \ d\Gamma
  }
  {
    \begin{aligned}
      a(u_h,v_h)= & \sum_{i=1}^2\Bigg(\int\limits_{\Omega_i} \alpha \nabla u_h \cdot \nabla v_h \ d\Omega_i \Bigg) \\
                  & - \int\limits_{\Gamma} \curlii[]{\alpha \nabla_{\bn} u_h} \llbracket v_h \rrbracket \ d\Gamma  \\
                  & - \int\limits_{\Gamma} \llbracket u_h \rrbracket \curlii[]{\alpha \nabla_{\bn} v_h} \ d\Gamma  \\
                  & + \int\limits_{\Gamma} \gamma \llbracket u_h \rrbracket \llbracket v_h \rrbracket \ d\Gamma
    \end{aligned}
  }
  \label{eq:N1_blinear}
\end{equation}
and
\begin{equation}
  F(v_h)=\sum_{i=1}^2\Bigg(\int\limits_{\Omega_i} f v_h \ d\Omega_i\Bigg),
  \label{eq:linear_form}
\end{equation}
where $\gamma\in \bR^+$ is the stabilization parameter.
In this formulation, the weighting parameters are defined as the measure fraction defined on element $K \in \pT_{h,\Gamma}$, e.g. $\beta_i = \text{meas}_d (K_i) / \text{meas}_d (K)$.
The formulation~\eqref{eq:NXFEM} can be shown to be consistent with the strong formulation~\eqref{eq:dd_poisson_eq} and also can be shown to be stable for a sufficiently large stabilization parameter.
The detailed analysis of this method is carried out in~\cite{hansbo_unfitted_2002}.
Under reasonable mesh assumptions, a priori error estimates for all $u \in H^1_0(\Omega)\cap H^2(\Omega_1 \cup \Omega_2)$ are given by
\begin{equation}
  \begin{aligned}
    \vertiii{u-u_h}_h       & \leqslant C \ h \sum_{i=1,2} \| u \|_{H^2(\Omega_i)},   \\
    \|u-u_h\|_{L^2(\Omega)} & \leqslant C \ h^2 \sum_{i=1,2} \| u \|_{H^2(\Omega_i)},
  \end{aligned}
  \label{eq:xfem_error}
\end{equation}
where the constant $C$ is completely independent of the location of the interface with respective to the mesh.
The mesh-dependent energy norms $\vertiii{\cdot}_h$ in the above estimates are defined as
\begin{equation}
  \whencolumns
  {
  \vertiii{v}^2_h := \| \nabla v\|^2_{L^2(\Omega_1 \cup \Omega_2)} + \|\llbracket v \rrbracket\|^2_{H^{\half}(\Gamma),h} + \left\| \curlii[]{ \nabla_{\bn} v}\right\|^2_{H^{-\half}(\Gamma),h},
  }
  {
  \begin{aligned}
    \vertiii{v}^2_h := & \| \nabla v\|^2_{L^2(\Omega_1 \cup \Omega_2)} + \|\llbracket v \rrbracket\|^2_{H^{\half}(\Gamma),h} \\
                       & + \left\| \curlii[]{ \nabla_{\bn} v}\right\|^2_{H^{-\half}(\Gamma),h},
  \end{aligned}
  }
  \label{eq:error_norms}
\end{equation}
and the mesh-dependent norms at the interface are given by
\whencolumns
{
\[
  \|v\|^2_{H^{\half}(\Gamma),h} = \sum_{K\in \pT_{h,\Gamma}} h^{-1} \left\| v\right\|^2_{L^2(\Gamma_K)} \quad \text{and} \quad
  \|v\|^2_{H^{-\half}(\Gamma),h} = \sum_{K\in \pT_{h,\Gamma}} h\left\|v \right\|^2_{L^2(\Gamma_K).}
\]

}
{
\[
  \|v\|^2_{H^{\half}(\Gamma),h} = \sum_{K\in \pT_{h,\Gamma}} h^{-1} \left\| v\right\|^2_{L^2(\Gamma_K)}
\]
and
\[
  \|v\|^2_{H^{-\half}(\Gamma),h} = \sum_{K\in \pT_{h,\Gamma}} h\left\|v \right\|^2_{L^2(\Gamma_K).}
\]
}
Thus, we have a formulation where the convergence of the method relies on the choice of the stabilization parameter.
Following the coercivity of the bilinear form \eqref{eq:N1_blinear}, as given in \ref{sect:appendix}, we have
\[
  \whencolumns
  {
    \begin{aligned}
      a(u,u) \geqslant & \sum_{K\in \widetilde{\pT}_{h}\setminus \pT_{h,\Gamma}} \|\alpha^{\half}\nabla u \|^2_{L^2(K)}
      + \sum_{K\in \pT_{h,\Gamma}} \half \|\alpha^{\half}\nabla u \|^2_{L^2(K)}
      + \frac{1}{\epsilon_c}\| \curlii[]{\alpha \nabla_{\bn}u} \|_{H^{-\half}(\Gamma),h }                                                     \\
                       & + \sum_{K\in\pT_{h,\Gamma}} \Big(\gamma - \frac{\epsilon_c}{h_K}\Big) \| \llbracket u \rrbracket\|^2_{L^2(\Gamma_K)}
      + \sum_{K\in\pT_{h,\Gamma}} \Big(\half - \frac{2C_\gamma}{\epsilon_c} \Big) \|\alpha^{\half} \nabla u\|_{L^2(K)}.
    \end{aligned}
  }
  {
    \begin{aligned}
      a(u,u) \geqslant & \sum_{K\in \widetilde{\pT}_{h}\setminus \pT_{h,\Gamma}} \|\alpha^{\half}\nabla u \|^2_{L^2(K)}                       \\
                       & + \sum_{K\in \pT_{h,\Gamma}} \half \|\alpha^{\half}\nabla u \|^2_{L^2(K)}                                            \\
                       & + \frac{1}{\epsilon_c}\| \curlii[]{\alpha \nabla_{\bn}u} \|_{H^{-\half}(\Gamma),h }                                  \\
                       & + \sum_{K\in\pT_{h,\Gamma}} \Big(\gamma - \frac{\epsilon_c}{h_K}\Big) \| \llbracket u \rrbracket\|^2_{L^2(\Gamma_K)} \\
                       & + \sum_{K\in\pT_{h,\Gamma}} \Big(\half - \frac{2C_\gamma}{\epsilon_c} \Big) \|\alpha^{\half} \nabla u\|_{L^2(K)}.
    \end{aligned}
  }
\]
This inequality utilizes Young's inequality for some ${\epsilon_c > 0}$ and follows trace inequality for all $K \in \pT_{h,\Gamma}$
\begin{equation}
  \| \curlii[]{\alpha \nabla_{\bn} u_{h}} \|^2_{L^2(\Gamma_K)} \leqslant \frac{C_{\gamma}}{h_K} \|\alpha^\half \nabla u_{h}\|^2_{L^2(K)}  .
  \label{eq:inverse_ineq}
\end{equation}
We can say that the bilinear form is coercive if the positivity of two terms in the last line in ensured, given by
$\epsilon_c \geqslant 4C_\gamma$, and $\gamma \geqslant \epsilon_c/h_K$.
Thus, the stabilization parameter can be given with the bound $\gamma \geqslant 4 C_\gamma/h_K$.

In the next section, we will discuss other modifications to the formulation (\ref{eq:NXFEM}), and approaches to estimate the stabilization parameters by estimating the constant $C_\gamma$.

\subsection{Variations in Nitsche's Formulations}
As mentioned in the previous section, Nitsche's formulation is stable only when the coercivity of the bilinear form is satisfied, i.e. by choosing a sufficiently large stabilization parameter $\gamma$.
Estimation of such stabilization parameter is a very delicate part of the method.
If the stabilization parameter is chosen to be too large, it gives rise to an ill-conditioned system matrix.
It becomes increasingly difficult to estimate the stabilization parameter for irregular meshes and higher order finite element discretization.

In addition to the stabilization parameter, the weighting parameters also play an important role in the stability of Nitsche's formulation.
A robust option for the weighting parameters is given in~\cite{barrau_robust_2012,annavarapu_robust_2012}, where it is suggested to use coefficient \(\alpha\) along with the measure of the cut-element.
The new weighting parameter is given as,
\begin{equation}
  \beta_i = \frac{\text{meas}_d(K_i)/\alpha_i } {\text{meas}_d(K_1)/\alpha_1 + \text{meas}_d(K_2)/\alpha_2}, \text{ for } i=1,2.
  \label{eq:barrau_annava_stab}
\end{equation}
This weighting parameter gives a better averaging for the discontinuous coefficients.

In this section, we explore different methods for estimating the stabilization parameter.

\subsubsection{Estimation of stabilization parameter - with Eigenvalue problem}
Here, we discuss the idea of estimating the stabilization parameter by solving a generalized eigenvalue problem.
This method was used for a particle-partition of unity method and later explored more for spline-based finite elements for harmonic and bi-harmonic problems~\cite{embar_imposing_2010,griebel_particle-partition_2003}.
In this method the value of the stabilization parameter is estimated by solving a generalized eigenvalue problem.
This approach is now widely used in the finite cell methods~\cite{schillinger_finite_2015,ruess_weakly_2013,ruess_weak_2014,jiang_robust_2015}, where Nitsche's method is used to enforce Dirichlet boundary condition.

We restate the weak formulation of Nitsche's method as, find $u_h \in \pV_h$, such that
\begin{equation}
  \label{eq:NXFEM-EV}
  a^N_1(u_h,v_h) = F(v_h), \qquad \forall v_h \in \pV_h,\tag{N-EV}
\end{equation}
where
\[
  \whencolumns
  {
    a^N_1(u_h,v_h) =  \sum_{i=1}^2\Bigg(
    \int\limits_{\Omega_i} \alpha \nabla u_h \cdot \nabla v_h \ d\Omega_i \Bigg)
    - \int\limits_{\Gamma} \curlii[]{\alpha \nabla_{\bn} u_h} \llbracket v_h \rrbracket \ d\Gamma
    - \int\limits_{\Gamma} \llbracket u_h \rrbracket \curlii[]{\alpha \nabla_{\bn} v_h} \ d\Gamma
    + \int\limits_{\Gamma} \gamma_1 \llbracket u_h \rrbracket \llbracket v_h \rrbracket \ d\Gamma.
  }
  {
    \begin{aligned}
      a^N_1(u_h,v_h) = & \sum_{i=1}^2\Bigg( \int\limits_{\Omega_i} \alpha \nabla u_h \cdot \nabla v_h \ d\Omega_i \Bigg) \\
                       & - \int\limits_{\Gamma} \curlii[]{\alpha \nabla_{\bn} u_h} \llbracket v_h \rrbracket \ d\Gamma   \\
                       & - \int\limits_{\Gamma} \llbracket u_h \rrbracket \curlii[]{\alpha \nabla_{\bn} v_h} \ d\Gamma   \\
                       & + \int\limits_{\Gamma} \gamma_1 \llbracket u_h \rrbracket \llbracket v_h \rrbracket \ d\Gamma.
    \end{aligned}
  }
  \label{eq:N2_blinear}
\]
Here, the weighting parameters $\beta_i$ are defined as in~\eqref{eq:barrau_annava_stab}.

The coercivity of the bilinear form relies on the trace inequality~\eqref{eq:inverse_ineq}.
A good estimate of $C_\gamma/h_K$ can be achieved by solving a generalized eigenvalue problem, as $C_\gamma/h_K$ is bounded from below by the largest eigenvalue of the auxiliary problem~\eqref{eq:eigval}.

We pose eigenvalue problems for each $K \in \pT_{h,\Gamma}$, and solve series of locally given element-wise problems, find $ \max(\lambda_K) \in \bR $ such that
\begin{equation}
  b_e(v_\lambda,v_\lambda) = \lambda_K c_e(v_\lambda,v_\lambda), \quad \forall v_\lambda \in \pV_h\vert_K
  \label{eq:eigval}
\end{equation}
where $\pV_h\vert_K$ is restriction of $\pV_h$ on a given element $K$.
Both bilinear forms on element $K \in \pT_{h,\Gamma}$ are defined as
\[
  b_e(v_\lambda,v_\lambda) = \int\limits_{\Gamma_K} \curlii[]{\alpha \nabla_{\bn} v_\lambda} \cdot \curlii[]{\alpha \nabla_{\bn} v_\lambda}\ d\Gamma,
\]
and
\[
  c_e(v_\lambda,v_\lambda) =  \sum_{i=1}^2 \int\limits_{K_i} \alpha \nabla v_\lambda \cdot \nabla v_\lambda \ d\Omega_i.
\]
In order to solve the generalized eigenvalue problem~\eqref{eq:eigval}, it is necessary that the bilinear form $c_e(\cdot,\cdot)$ has only trivial kernel.
This can be achieved, if the function space $\pV_h\vert_K$ is defined in the space of polynomials which are orthogonal to constants.
From the construction, it is clear that the bilinear form $c_e(\cdot,\cdot)$, is a representation of a local stiffness matrix and here the kernel of the local stiffness matrix is known to be a constant vector.
Algebraically, we can use a deflation method to eliminate the influence of the trivial kernel from matrix representation of both bilinear forms, $b_e(\cdot,\cdot)$ and $c_e(\cdot,\cdot)$, and still retain other spectral properties of matrices.
Thus, solving the generalized eigenvalue problem of the deflated system is equivalent to solving the original eigenvalue problem~\eqref{eq:eigval}, and we can use the largest eigenvalue in the estimation of the stabilization parameter.
To ensure the boundedness of~\eqref{eq:inverse_ineq}, we take the value of element-wise stabilization parameter 4 times larger than the largest eigenvalue.
Hence, the stabilization parameter is computed element-wise as, $\gamma_1 = 4\max(\lambda_K)$ to satisfy the condition $\gamma_1 \geqslant 4C_\gamma/h_K$.

Above, we have defined local eigenvalue problems for each element $K \in \pT_{h,\Gamma}$.
An alternate option is to create a global system for all the cut-elements and solving a global eigenvalue problem.
In this case, the stabilization parameter is estimated by the largest eigenvalue of the global system.
For irregular meshes and complex domains, an interface can intersect the mesh arbitrarily and a very small cut in one element can influence the largest eigenvalue in the global setting.
Hence, we choose to solve a series of local eigenvalue problems and estimate the local stabilization parameter for each cut-element.
This approach is more beneficial, as firstly we avoid computing the largest eigenvalue of the global problem and secondly the effect of small cut-elements is localized.

\subsubsection{Estimation of stabilization parameter - with Lifting operator}
To avoid the estimation process of the stabilization parameter for Nitsche's formulation, an alternative method is proposed in~\cite{lehrenfeld_removing_2016,lehrenfeld_space-time_2015}.
In this method, the stabilization parameter is chosen locally in implicit manner, similar to the previous section.
This strategy is common in the DG method~\cite{bassi_high-order_1997}.
The stability and error analysis of the DG discretization equipped with the lifting operators is carried out in~\cite{brezzi_discontinuous_2000}.

First, we introduce an element-wise lifting operator, $\pazocal{L}_K(\cdot)$, which lifts the functions defined on the cut-elements into the space of polynomials which are orthogonal to constants, $\pazocal{L}_K:\pV_h\vert_K \to \pazocal{W}_h$, where
\begin{equation*}
  \pazocal{W}_h := \{u_h\in L^2(K): u_h|_{K} \in \pP_1 \cap (\pP_0)^\perp,\ \forall K \in \pT_{h,\Gamma}\}.
\end{equation*}
On the uncut elements, i.e. $K \in \widetilde{\pT}_{h} \setminus \pT_{h,\Gamma}$ the lifting operator is defined as $\pazocal{L}_K(u_h) = 0,$.
While on the cut-elements, the lifting operator is defined as, find $w_h:=\pazocal{L}_K(u_h) \in \pazocal{W}_h$ such that
\begin{equation*}
  b_l(w_h,v_h) = c_l(u_h,v_h), \quad \forall u_h,v_h \in \pV_h\vert_K.
\end{equation*}
The bilinear forms on an element $K \in \pT_{h,\Gamma}$ are defined as
\[
  b_l(w_h,v_h) =  \sum_{i=1}^2 \int\limits_{K_i} \alpha \nabla w_h \cdot \nabla v_h \ d\Omega_i
\]
and
\[
  c_l(u_h,v_h) = - \int\limits_{\Gamma_K} \llbracket u_h \rrbracket \curlii[]{\alpha \nabla_{\bn} v_h} \ d\Gamma.
\]
The coercivity of the bilinear form of the original formulation~\eqref{eq:NXFEM} can be ensured if an additional term stemming from the lifting operators is added
\begin{equation*}
  \whencolumns
  {
    \begin{aligned}
      \sum_{i=1}^2\Bigg(
      \int\limits_{\Omega_i} \alpha \nabla u_h \cdot \nabla u_h \ d\Omega_i \Bigg)
      - 2 \int\limits_{\Gamma} \llbracket u_h \rrbracket \curlii[]{\alpha \nabla_{\bn} u_h} \ d\Gamma
       & + 2 \sum_{i=1}^2 \sum_{K \in \pT^i_{h,\Gamma} } \int\limits_{K_i} \alpha \pazocal{L}_K (\nabla u_h) \cdot \pazocal{L}_K(\nabla u_h)\ d\Omega_i \\
       & \geqslant \half
      \sum_{i=1}^2 \Bigg(\int\limits_{\Omega_i} \alpha \nabla u_h \cdot \nabla u_h \ d\Omega_i \Bigg). \end{aligned}
  }
  {
    \begin{aligned}
      \sum_{i=1}^2\Bigg(
       & \int\limits_{\Omega_i} \alpha \nabla u_h \cdot \nabla u_h \ d\Omega_i \Bigg)
      - 2 \int\limits_{\Gamma} \llbracket u_h \rrbracket \curlii[]{\alpha \nabla_{\bn} u_h} \ d\Gamma                                                   \\
       & + 2 \sum_{i=1}^2 \sum_{K \in \pT^i_{h,\Gamma} } \int\limits_{K_i} \alpha \pazocal{L}_K (\nabla u_h) \cdot \pazocal{L}_K(\nabla u_h)\ d\Omega_i \\
       & \geqslant \half
      \sum_{i=1}^2 \Bigg(\int\limits_{\Omega_i} \alpha \nabla u_h \cdot \nabla u_h \ d\Omega_i \Bigg). \end{aligned}
  }
\end{equation*}
Addition of such a term in the bilinear form ensures the coercivity for any positive stabilization parameter.
The updated weak formulation can be given as, find $u_h \in \pV_h$, such that
\begin{equation}
  \label{eq:NXFEM-LO}
  a^N_2(u_h,v_h) = F(v_h), \qquad \forall v_h \in \pV_h,\tag{N-LO}
\end{equation}
where
\[
\whencolumns
  {
    \begin{aligned}
      a^N_2(u_h,v_h) = \sum_{i=1}^2\Bigg(
      \int\limits_{\Omega_i} \alpha \nabla u_h \cdot \nabla v_h \ d\Omega_i \Bigg)
       & - \int\limits_{\Gamma} \curlii[]{\alpha \nabla_{\bn} u_h} \llbracket v_h \rrbracket \ d\Gamma
      - \int\limits_{\Gamma} \llbracket u_h \rrbracket \curlii[]{\alpha \nabla_{\bn} v_h} \ d\Gamma + \int\limits_{\Gamma} \gamma_2 \llbracket u_h \rrbracket \llbracket v_h \rrbracket \ d\Gamma \\
       & + 2 \sum_{i=1}^2 \sum_{K \in \pT^i_{h,\Gamma} } \int\limits_{K_i} \alpha \pazocal{L}_K(\nabla u_h) \cdot  \pazocal{L}_K(\nabla v_h)\ d\Omega_i.
    \end{aligned}
  }
  {
    \begin{aligned}
      a^N_2(u_h,v_h) & = \sum_{i=1}^2\Bigg(
      \int\limits_{\Omega_i} \alpha \nabla u_h \cdot \nabla v_h \ d\Omega_i \Bigg)                                                                                                                                 \\
                     & - \int\limits_{\Gamma} \curlii[]{\alpha \nabla_{\bn} u_h} \llbracket v_h \rrbracket \ d\Gamma                                                                                               \\
                     & - \int\limits_{\Gamma} \llbracket u_h \rrbracket \curlii[]{\alpha \nabla_{\bn} v_h} \ d\Gamma + \int\limits_{\Gamma} \gamma_2 \llbracket u_h \rrbracket \llbracket v_h \rrbracket \ d\Gamma \\
                     & + 2 \sum_{i=1}^2 \sum_{K \in \pT^i_{h,\Gamma} } \int\limits_{K_i} \alpha \pazocal{L}_K(\nabla u_h) \cdot  \pazocal{L}_K(\nabla v_h)\ d\Omega_i.                                             \\
    \end{aligned}
  }
  \label{eq:N4_blinear}
\]
In the original work~\cite{lehrenfeld_removing_2016} different weighting parameters and the stabilization parameter are chosen.
In order to improve the robustness of the formulation for the problems with highly varying coefficients, we choose the weighting parameters defined by~\eqref{eq:barrau_annava_stab} and the stabilization parameter is chosen as in~\cite{annavarapu_robust_2012}, hence
\[
  \gamma_2 = \frac{\text{meas}_{d-1}(\Gamma_K)}{\text{meas}_d(K_1)/\alpha_1 + \text{meas}_d(K_2)/\alpha_2} \ .
\]

\subsubsection{Ghost Penalty Stabilization}
Ghost penalty was introduced for the fictitious domain methods as an additional stabilization term to enhance the robustness of Nitsche's method irrespective of the interface location~\cite{burman_ghost_2010}.
The idea of such a stabilization term was used for the problems with dominant transport to penalize the jumps in the normal derivative across the interior faces of elements~\cite{douglas_interior_1976}.
This method was recently applied in the context of convection - diffusion - reaction problem and Stoke's problem~\cite{burman_edge_2004,Burman_edge_2006}.
This kind of stabilization term is also applied in the XFEM context for incompressible elasticity problems to penalize the  jump in pressure~\cite{becker_nitsche_2009}.
In the fictitious domain methods, it was shown that under certain conditions, especially for sliver cuts, the error bounds of Nitsche's formulation degrade and it may even produce diverging solutions.
But if the ghost penalty stabilization term is added to the formulation, Nitsche's method becomes stable even for pathological cases~\cite{de_prenter_note_2018}.

For the standard Nitsche formulation, as we have mentioned earlier, the condition number of the system matrix depends on the cut position and the conditioning of the system can become arbitrarily bad when an interface passes very close to element faces or nodes.
The ghost penalty method overcomes this issue by extending the coercivity of the bilinear form, from the computational domain to the fictitious part.
By adding this term, the condition number of the system can be bounded independent of location of the interface in the mesh.

We define the set of faces $\pazocal{G}$ for each subdomain $\Omega_i$,
\[
  \pazocal{G}_{\Gamma,i} =\{ G \subset \partial K : K\in \pT^i_{h,\Gamma},\ \partial K \cap \partial \pT_i = \emptyset\}, \quad i=1,2.
\]
This set includes all the faces which are associated with the cut-elements, except the ones on boundary.
The ghost penalty term is defined as,
\begin{equation}
  \begin{aligned}
    g(u_h,v_h) = \sum_{i=1}^2 \sum_{G\in \pazocal{G}_{\Gamma,i}} \int\limits_{G} \epsilon_G h_G \alpha \llbracket \nabla_{\bn} u_h \rrbracket \llbracket \nabla_{\bn} v_h \rrbracket \ dG,
  \end{aligned}
\end{equation}
where $h_G$ is the diameter of face $G$, and $\epsilon_G$ is a positive constant.

In the previous section, we observed that the choice of the averaging function is also quite important.
For the problems with highly varying coefficients this stabilization term is an attractive option, when the averaging function does not provide sufficient stability.
In addition, for highly varying coefficients~\cite{burman_numerical_2012}, the definition of the $\beta_i$, is changed to
\begin{equation}
  \beta_1=\frac{\alpha_2}{\alpha_1+\alpha_2} \quad \text{and} \quad \beta_2 = \frac{\alpha_1}{\alpha_1+\alpha_2},
\end{equation}
and the stabilization parameter can also be chosen to be coefficient dependent
\[
  \gamma_3 = \gamma_0\frac{2 \alpha_1 \alpha_2}{ \alpha_1+\alpha_2}, \quad \text{where } \gamma_0 \in \bR^+.
\]
Now, the weak formulation equipped with the ghost penalty term is given as, find $u_h \in \pV_h$, such that
\begin{equation}
  \label{eq:NXFEM-GP}
  a^N_3(u_h,v_h) = F(v_h), \qquad \forall v_h \in \pV_h,\tag{N-GP}
\end{equation}
where
\[
  \whencolumns
  {
    \begin{aligned}
      a^N_3(u_h,v_h) = \sum_{i=1}^2\Bigg(\int\limits_{\Omega_i}  \alpha \nabla u_h \cdot \nabla v_h \ d\Omega_i \Bigg)
       & - \int\limits_{\Gamma} \curlii[]{\alpha \nabla_{\bn} u_h} \llbracket v_h \rrbracket \ d\Gamma
      - \int\limits_{\Gamma} \llbracket u_h \rrbracket \curlii[]{\alpha \nabla_{\bn} v_h} \ d\Gamma
      + \int\limits_{\Gamma} \frac{\gamma_3}{h} \llbracket u_h \rrbracket \llbracket v_h \rrbracket \ d\Gamma                                                                         \\
       & + \sum_{i=1}^2 \sum_{G \in \pazocal{G}_{\Gamma,i}} \int\limits_{G} \epsilon_G h_G \alpha \llbracket \nabla_{\bn} u_h \rrbracket \llbracket \nabla_{\bn} v_h \rrbracket \ dG.
    \end{aligned}
  }
  {
    \begin{aligned}
      a^N_3(u_h,v_h) & = \sum_{i=1}^2\Bigg(
      \int\limits_{\Omega_i}  \alpha \nabla u_h \cdot \nabla v_h \ d\Omega_i \Bigg)                                                                                                                 \\
                     & - \int\limits_{\Gamma} \curlii[]{\alpha \nabla_{\bn} u_h} \llbracket v_h \rrbracket \ d\Gamma                                                                                \\
                     & - \int\limits_{\Gamma} \llbracket u_h \rrbracket \curlii[]{\alpha \nabla_{\bn} v_h} \ d\Gamma
      + \int\limits_{\Gamma} \frac{\gamma_3}{h} \llbracket u_h \rrbracket \llbracket v_h \rrbracket \ d\Gamma                                                                                       \\
                     & + \sum_{i=1}^2 \sum_{G \in \pazocal{G}_{\Gamma,i}} \int\limits_{G} \epsilon_G h_G \alpha \llbracket \nabla_{\bn} u_h \rrbracket \llbracket \nabla_{\bn} v_h \rrbracket \ dG.
    \end{aligned}
  }
  \label{eq:N5_blinear}
\]
The error analysis of this method was carried out in~\cite{burman_numerical_2012} and it was shown to have optimal convergence rates in the $L^2$-norm and the mesh-dependent energy norm.

\begin{remark}
  In the numerical experiments, we choose the value of the stabilization parameter \(\gamma_0 = 10\), and value of the constant in the ghost penalty term as \(\epsilon_G = 0.1\).
\end{remark}

 \section{Multigrid for Nitsche-XFEM} \label{sect:multigrid}
As mentioned in the introduction, a multigrid method is a combination of smoothing iterations and coarse level corrections.
Coarse level corrections are heavily dependent on well chosen transfer operators which can be used to restrict residual from a fine to a coarse level and prolongate correction from a coarse to a fine level.
In this section, we introduce our multigrid method for the XFEM discretization and discuss a new transfer operator based on the $L^2$-projection and pseudo-$L^2$-projection in the XFEM framework.
We note, this multilevel method relies only on the unfitted meshes and the enriched XFEM spaces, it is agnostic of any method chosen to enforce the boundary or the interface conditions.

\subsection{Multilevel Framework for XFEM discretization}
In this subsection, we provide a framework for creating a multilevel XFEM space hierarchy from the hierarchy of background meshes.

Let us define a mesh hierarchy of background meshes for levels, $\ell$, where $\ell \in \{0,\dots, L\}$.
The coarsest level is denoted as $\ell=0$ and the finest level is denoted by $\ell=L$.
The original background mesh is denoted as the mesh on the finest level, ${\widetilde{\pT}_L := \widetilde{\pT}_h}$.
We define the coarse level mesh hierarchy as, ${\widetilde{\pT}_{\ell}}$, ${\ell \in \{0,\dots,L-1\}}$.
At each level $\ell$, we require that the domain $\Omega$ is encapsulated by the mesh hierarchy, ${\Omega \subset (\widetilde{\pT}_\ell)_{\ell=\{0,\dots,L-1\}}}$.
Now, on the mesh hierarchy, we define a finite element space associated with each of these background mesh,
\[
  \whencolumns{
    \widetilde{\pV}_\ell = \{v\in H^1(\Omega) : v\vert_{K} \in \pP_1,\ v\vert_{\partial \Omega}=0,\ \forall K\in \widetilde{\pT}_\ell \}, \quad \forall  \ell\in \{0,\dots,L-1\}.
  }
  {
    \begin{aligned}
      \widetilde{\pV}_\ell = \{v\in H^1(\Omega) : v\vert_{K} \in \pP_1,\ v\vert_{\partial \Omega}=0,\ \forall K\in \widetilde{\pT}_\ell \}, \\
      \forall  \ell\in \{0,\dots,L-1\}.
    \end{aligned}
  }
\]
On the finest level the definition of the finite element space is taken directly from the original problem, i.e. ${\widetilde{\pV}_L := \widetilde{\pV}_h}$.
When the hierarchy of the meshes is nested, the hierarchy of the FE spaces associated with these meshes are also nested,
\[
  \widetilde{\pV}_{\ell-1} \subset \widetilde{\pV}_\ell, \quad \forall  \ell \in \{1,\dots,L\}.
\]

In the XFEM discretization, the background mesh is enriched, decomposed and then it is associated with each subdomain.
In the multigrid framework this procedure is carried out on each level in the mesh hierarchy,
\[
  \whencolumns
  {
    \pT_{\ell,i} = \{K \in \widetilde{\pT}_\ell : K \cap \Omega_i \neq \emptyset\}, \quad \forall  \ell\in \{0,\dots,L-1\}, \ i\in\{1,2\}.
  }
  {
    \begin{aligned}
      \pT_{\ell,i} = \{K \in \widetilde{\pT}_\ell: & K \cap \Omega_i \neq \emptyset\},                \\
                                                   & \forall  \ell\in \{0,\dots,L-1\}, \ i\in\{1,2\}.
    \end{aligned}
  }
\]
The mesh hierarchy can be achieved by uniform coarsening, if possible, of the background mesh of the original discretization.

\begin{figure*}[t]
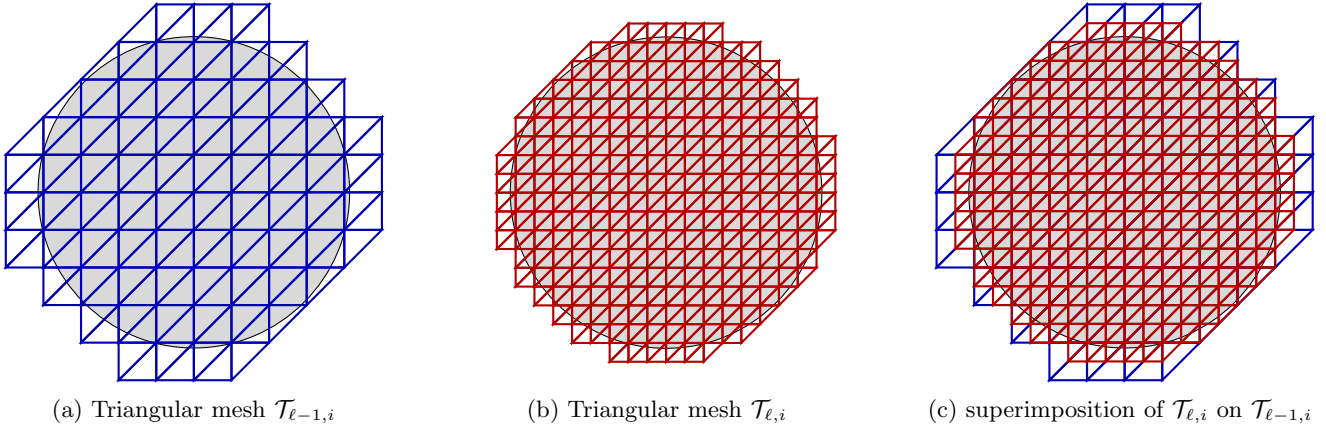

  \centering
  \subcaptionbox{Triangular mesh ${\pT}_{\ell-1,i}$}[.3\textwidth]{

  }
  \hfill
  \subcaptionbox{Triangular mesh $\pT_{\ell,i}$}[.3\textwidth]{%
  }
  \hfill
  \subcaptionbox{superimposition of $\pT_{\ell,i}$ on $\pT_{\ell-1,i}$}[.3\textwidth]{
    %
  }
  \caption[2D mesh]{2D Triangular meshes on different levels encapsulating the domain $\Omega_i$, (domain $\Omega_i$ is shaded in gray)}
  \label{fig:multilevel_decomp}
\end{figure*}

Now, we exploit the definition of the Heaviside function (\ref{eq:Heaviside}) to restrict the finite element space $(\widetilde{\pV}_\ell)_{\ell\in\{0,\dots,L-1\}}$ in each subdomain $\Omega_i$, as
\[
  \pV_{\ell,i} = \mathcal{H}_i(\bx)\widetilde{\pV}_\ell, \quad \forall  \ell \in \{0,\dots,L-1\},\ i \in \{1,2\}.
\]
Similar to the finest level mesh, we borrow the definition of FE space from the finest level on each subdomain, i.e. ${\pV_{L,i} := \pV_{h,i}}$.
In the XFEM framework, even if we have a nested mesh hierarchy of the background mesh, due to the arbitrary location of interfaces, we could have a non-nested mesh hierarchy for each subdomain.
We can see an example for this in Figure~\ref{fig:multilevel_decomp}.
There, the meshes associated with the domain $\Omega_i$ at different levels are not nested.
As the meshes are not nested, the finite element spaces associated with the meshes are also not nested, ${\pV_{\ell-1,i} \not\subset \pV_{\ell,i}}$.

\begin{figure*}[t]
  \centering
  \begin{tikzpicture}[xscale=1.0,yscale=1.6]

\begin{scope}[shift={(-0.5,0)}]
  \foreach \i in {4,...,11}{
    \draw[gray,line width=0.7mm](\i,-0.05) -- (\i,0.05);
  }
  \draw[gray,line width=0.7mm](4,0) -- (11.1556,0);
  \draw[](3,0.5)--(3,0.5) node {$L$};
  
  \draw[line width=0.7mm] (4,1)--(4,1) node [left = 0.01 cm] {$1$};
  \draw[line width=0.7mm] (4,0)--(4,0) node [left = 0.01 cm] {$0$};

  \draw[blue!70!black,line width=0.7mm](4,1) -- (5,0);
  \draw[blue!70!black,line width=0.7mm](4,0) -- (5,1) -- (6,0);
  \draw[blue!70!black,line width=0.7mm](5,0) -- (6,1) -- (7,0);
  \draw[blue!70!black,line width=0.7mm](6,0) -- (7,1) -- (8,0);
  \draw[blue!70!black,line width=0.7mm](7,0) -- (8,1) -- (9,0);
  \draw[blue!70!black,line width=0.7mm](8,0) -- (9,1) -- (10,0);
  \draw[blue!70!black,line width=0.7mm](9,0) -- (10,1) -- (11,0);
  \draw[blue!70!black,line width=0.7mm](10,0) -- (11,1) -- (11.1556,0.8444);
  \draw[blue!70!black,line width=0.7mm](11,0) -- (11.1556,0.1556);

  \draw[line width=0.7mm] (11.1556,-0.2)--(11.1556,1.2) node [above = 0.01 cm] {$\Gamma$};
\end{scope}

\begin{scope}[shift={(0.5,0)}]
  \foreach \i in {12,...,16}{
    \draw[gray,line width=0.7mm](\i,-0.05) -- (\i,0.05);
  }
  \draw[gray,line width=0.7mm](11.1556,0) -- (16,0);

  \draw[blue!70!black,line width=0.7mm](11.1556,0.8444) -- (12,0);
  \draw[blue!70!black,line width=0.7mm](11.1556,0.1556) -- (12,1) -- (13,0);
  \draw[blue!70!black,line width=0.7mm](12,0) -- (13,1) -- (14,0);
  \draw[blue!70!black,line width=0.7mm](13,0) -- (14,1) -- (15,0);
  \draw[blue!70!black,line width=0.7mm](14,0) -- (15,1) -- (16,0);
  \draw[blue!70!black,line width=0.7mm](15,0) -- (16,1);
  
  \draw[line width=0.7mm] (11.1556,-0.2)--(11.1556,1.2) node [above = 0.01 cm] {$\Gamma$} ;
\end{scope}

\begin{scope}[shift={(-0.5,-1.5)}]
  \foreach \i in {4,...,11}{
    \draw[gray,line width=0.7mm](\i,-0.05) -- (\i,0.05);
  }
  \draw[](3,0.5)--(3,0.5) node {$L-1$};
  
  \draw[line width=0.7mm,gray](4,0) -- (11.1556,0);
  \draw[line width=0.7mm] (4,1)--(4,1) node [left = 0.01 cm] {$1$};
  \draw[line width=0.7mm] (4,0)--(4,0) node [left = 0.01 cm] {$0$};
 
\draw[blue!70!black,line width=0.7mm](4,1) -- (6,0);
  \draw[blue!70!black,line width=0.7mm](4,0) -- (6,1) -- (8,0);
  \draw[blue!70!black,line width=0.7mm](6,0) -- (8,1) -- (10,0);
  \draw[blue!70!black,line width=0.7mm](8,0) -- (10,1) -- (11.1556,0.4222);
  \draw[blue!70!black,line width=0.7mm](10,0) -- (11.1556,0.5778);

  \draw[red!70!black,line width=0.7mm](10,1) -- (11,0.2161) -- (11.1556,0.2482);
  \draw[red!70!black,line width=0.7mm](10,0) -- (11,0.7839) -- (11.1556,0.7518);

  \draw[line width=0.7mm] (11.1556,1.2)--(11.1556,-0.2);
\end{scope}

\begin{scope}[shift={(0.5,-1.5)}]
  \foreach \i in {12,...,16}{
    \draw[gray,line width=0.7mm](\i,-0.05) -- (\i,0.05);
  }
  \draw[line width=0.7mm,gray](11.1556,0) -- (16,0);

  \draw[line width=0.7mm,red!70!black](11.1556,0.3565) -- (12,0);
  \draw[line width=0.7mm,red!70!black](11.1556,0.6435) -- (12,1);

  \draw[blue!70!black,line width=0.7mm](11.1556,0.4222) -- (12,0);
  \draw[blue!70!black,line width=0.7mm](11.1556,0.5778) -- (12,1) -- (14,0);
  \draw[blue!70!black,line width=0.7mm](12,0) -- (14,1) -- (16,0);
  \draw[blue!70!black,line width=0.7mm](14,0) -- (16,1);

  \draw[line width=0.7mm] (11.1556,1.2)--(11.1556,-0.2);
\end{scope}

\begin{scope}[shift={(-0.5,-3)}]
  \foreach \i in {4,...,11}{
    \draw[gray,line width=0.7mm](\i,-0.05) -- (\i,0.05);
  }
  \draw[](3,0.5)--(3,0.5) node {$L-2$};

  \draw[line width=0.7mm] (4,1)--(4,1) node [left = 0.01 cm] {$1$};
  \draw[line width=0.7mm] (4,0)--(4,0) node [left = 0.01 cm] {$0$};

  \draw[gray,line width=0.7mm](4,0) -- (11.1556,0);
\draw[blue!70!black,line width=0.7mm](4,1) -- (8,0);
  \draw[blue!70!black,line width=0.7mm](4,0) -- (8,1) -- (11.1556,0.2111);
  \draw[blue!70!black,line width=0.7mm](8,0) -- (11.1556,0.7889);

  \draw[red!70!black,line width=0.7mm](8,1) -- (9,0.7133) -- (10,0.4266) -- (11,0.2577)--  (11.1556,0.2646);
  \draw[red!70!black,line width=0.7mm](8,0) -- (9,0.2867) -- (10,0.5734)--(11,0.7423) -- (11.1556,0.7354);
  
  \draw[line width=0.7mm] (11.1556,1.2)--(11.1556,-0.2);
\end{scope}

\begin{scope}[shift={(0.5,-3)}]
  \foreach \i in {12,...,16}{
    \draw[gray,line width=0.7mm](\i,-0.05) -- (\i,0.05);
  }
  \draw[gray,line width=0.7mm](11.1556,0) -- (16,0);

  \draw[line width=0.7mm,red!70!black](11.1556,0.0753) -- (12,0);
  \draw[line width=0.7mm,red!70!black](11.1556,0.9247) -- (12,1);
 
  \draw[line width=0.7mm,blue!70!black](11.1556,0.2111) -- (12,0);
  \draw[line width=0.7mm,blue!70!black](11.1556,0.7889) -- (12,1) -- (16,0);
  \draw[line width=0.7mm,blue!70!black](12,0) -- (16,1);

  \draw[line width=0.7mm] (11.1556,1.2)--(11.1556,-0.2);
\end{scope}

   \end{tikzpicture}
  \caption[nested basis]{Basis of the cut domain in $1D$, in blue we see the original basis $\pV_{L,i}\not\supset\pV_{L-1,i}\not\supset\pV_{L-2,i}$, in red we see the basis computed with $L^2$-projections which are nested, $\pX_{L,i}\supset\pX_{L-1,i}\supset\pX_{L-2,i}$}
  \label{fig:basis_l2}
\end{figure*}
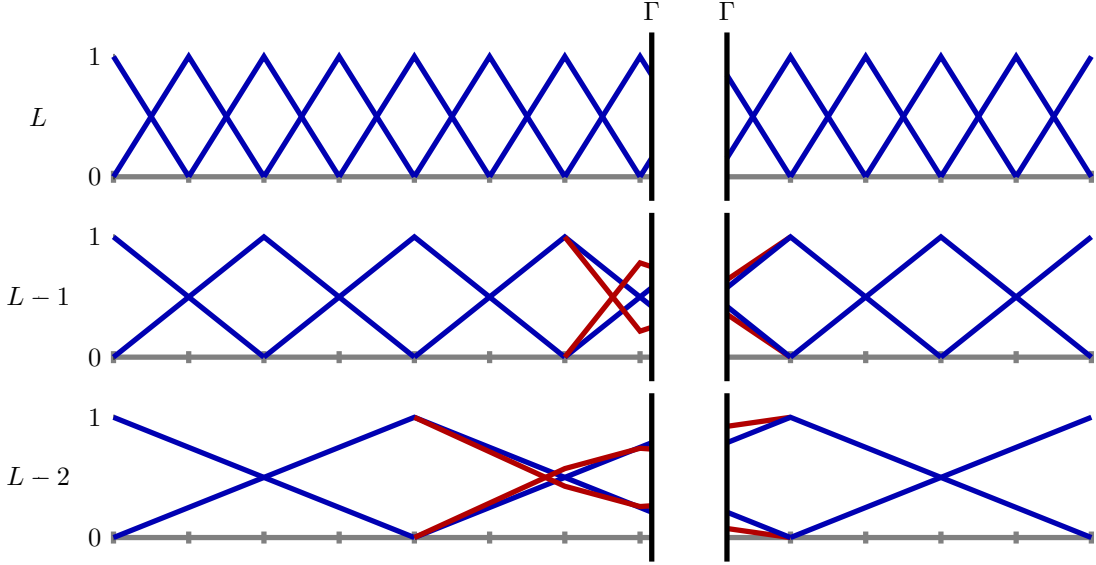

In this work, we want to create the hierarchy of nested FE spaces associated with each subdomain.
We define the prolongation operators which connects the spaces, ${(\pV_{\ell,i})_{\ell \in \{0,\ldots,L\}}}$, as
\begin{equation}
  \Pi_{\ell-1,i}^\ell : \pV_{\ell-1,i} \to \pV_{\ell,i}, \ \ \forall  \ell \in\{1,\dots,L\}, i\in\{1,2\},
  \label{eq:prolongation}
\end{equation}
such that ${\Pi_{\ell-1,i}^{\ell} \pV_{\ell-1,i} \subset \pV_{\ell,i}}$.
Now using this prolongation operator, we can construct a FE space associated with mesh $\pT_{\ell,i}$ by composition of the prolongation operators,
\[
  \pX_{\ell,i} := \Pi_{L-1,i}^L \cdots \Pi_{\ell,i}^{\ell+1} \pV_{\ell,i}, \ \forall  \ell \in \{1,\dots,L-1\}, i\in\{1,2\}.
\]
On the finest level we define ${\pX_{L,i} := \pV_{L,i}}$, as the definition of FE space on the finest level is kept untouched.
By recursive application of the prolongation operator we can create the hierarchy of the nested spaces as
\[
  \whencolumns{
    \underbrace { (\Pi_{L-1,i}^L \cdots \Pi_{\ell,i}^{\ell+1})(\Pi_{\ell-1,i}^\ell \pV_{\ell-1,i})}_{=:\pX_{\ell-1,i}}
    \subset
    \underbrace { (\Pi_{L-1,i}^L \cdots \Pi_{\ell,i}^{\ell+1})(\pV_{\ell,i})}_{=:\pX_{\ell,i}}, \quad \forall  \ell \in \{ 1,\dots,L-1\},\ i\in\{1,2\}.
  }{
    \begin{aligned}
      \underbrace { (\Pi_{L-1,i}^L \cdots \Pi_{\ell,i}^{\ell+1})(\Pi_{\ell-1,i}^\ell \pV_{\ell-1,i})}_{=:\pX_{\ell-1,i}}
      \subset
      \underbrace { (\Pi_{L-1,i}^L \cdots \Pi_{\ell,i}^{\ell+1})(\pV_{\ell,i})}_{=:\pX_{\ell,i}}, \\ \quad \forall  \ell \in \{ 1,\dots,L-1\},\ i\in\{1,2\}.
    \end{aligned}}
\]
Thus, the composition of the prolongation operators applied on the finest level FE space generates a nested space hierarchy,
\[
  \whencolumns{
    \pX_{0,i} \subset \pX_{1,i} \subset \cdots \subset \pX_{\ell,i}\subset \cdots \subset \pX_{L-1,i} \subset \pX_{L,i}, \quad  i\in\{1,2\}.
  }
  {
    \begin{aligned}
      \pX_{0,i} \subset \pX_{1,i} \subset \cdots \subset \pX_{\ell,i}\subset \cdots \subset \pX_{L-1,i} \subset \pX_{L,i},  i\in\{1,2\}.
    \end{aligned}
  }
\]
For simplicity and compactness, we define the prolongation operator for the whole domain as
\begin{equation}
  \Pi_{\ell-1}^\ell := \Pi_{\ell-1,1}^\ell \oplus \Pi_{\ell-1,2}^\ell.
\end{equation}
This prolongation operator inherits the same properties from its counterparts defined on each subdomain.
A hierarchy of nested spaces for the whole domain can be generated with the same procedure on the enriched FE space, for all ${\ell \in\{0,\dots,L-1\}}$, i.e.
\begin{equation*}
  \whencolumns{
    \pX_{\ell}  = \Pi_{\ell-1}^\ell \underbrace{ (\pV_{\ell,1}\oplus \pV_{\ell,2})}_{=:\pV_\ell} = \underbrace{ (\Pi_{L-1,1}^L \cdots \Pi_{\ell,1}^{\ell+1} \pV_{\ell,1}) \oplus( \Pi_{L-1,2}^L \cdots \Pi_{\ell,2}^{\ell+1} \pV_{\ell,2})}_{=\pX_{\ell,1} \oplus \pX_{\ell,2}}.
  }
  {
    \begin{aligned}
      \pX_{\ell} & = \Pi_{\ell-1}^\ell \underbrace{ (\pV_{\ell,1}\oplus \pV_{\ell,2})}_{=:\pV_\ell} \\ &= \underbrace{ (\Pi_{L-1,1}^L \cdots \Pi_{\ell,1}^{\ell+1} \pV_{\ell,1}) \oplus( \Pi_{L-1,2}^L \cdots \Pi_{\ell,2}^{\ell+1} \pV_{\ell,2})}_{=\pX_{\ell,1} \oplus \pX_{\ell,2}}.
    \end{aligned}
  }
\end{equation*}
\whencolumns{
  Thus, we can create the sequence of nested spaces $(\pX_{\ell})_{\ell\in{0,\ldots,L}}$.
}
{
  Thus, we can create the sequence of nested spaces\linebreak $(\pX_{\ell})_{\ell\in{0,\ldots,L}}$.
}
In Figure~\ref{fig:basis_l2}, we can see how the basis functions created by the nested FE spaces differ from the non-nested FE spaces.

\subsection{Variational Transfer}
In this section, we discuss the computation of the prolongation operator for the XFEM framework.
In the context of non-conforming domain decomposition methods and contact problems, the information transfer between non-conforming meshes is realized by means of global $L^2$-projections~\cite{belgacem_mortar_1999,wohlmuth_mortar_2000}.
The mortar methods were introduced to couple different discretization on the interfaces of subdomains, where the meshes of the subdomains does not necessarily match at the interface.
In the mortar method the $L^2$-projections are performed to couple these discretizations on the trace spaces defined on the boundaries.
We exploit the same strategy to couple domains between different FE spaces, i.e. we create the transfer operators between the successive meshes in the multilevel hierarchy.
In the theory of mortar method, a mortar side and a non-mortar side is chosen, the projection operator maps a function from a mortar side to a non-mortar side.
Thus, a mortar side and a non-mortar side can be regarded as a domain and an image of a projection operator.
In the multigrid framework, we associate a fine space with the non-mortar side, while a coarse space with the mortar side, as we aim to compute the prolongation operator from a coarse to a fine space.
The projection operator between spaces \( (\pV_\ell)_{\ell \in \{0,\ldots,L\}} \) in our mesh hierarchy is defined as in \eqref{eq:prolongation}.
For example, the operator ${\Pi_{\ell-1}^\ell}$ projects the element ${v \in \pV_{\ell-1}}$ into $\pV_\ell$, given as ${\Pi_{\ell-1}^\ell v:= w}$.

Setting ${(I=\pT_{\ell-1,i} \cap \pT_{\ell,i})_{i=1,2}}$, we define our projection operator, \(\Pi_{\ell-1}^{\ell}: \pV_{\ell-1} \to \pV_{\ell}\) as
\begin{equation}
  \Pi_{\ell-1}^{\ell}v\in \pV_\ell:\ (\Pi_{\ell-1}^\ell v, \mu)_{L^2(I)} = ( w, \mu )_{L^2(I)},  \quad \forall \mu \in \pM_\ell,
  \label{eq:variational_coupling}
\end{equation}
where $\pM_\ell$ is a space of Lagrange multiplier defined on the fine level, thus we also have $\dim(\pM_\ell) =\dim(\pV_\ell)$.
Reformulating \eqref{eq:variational_coupling}, we get weak equality condition on the intersection of meshes
\begin{equation}
  (v-\Pi_{\ell-1}^\ell v , \mu )_{L^2(I)} = (v-w ,\mu)_{L^2(I)} = 0,\quad \forall \mu \in \pM_\ell.
  \label{eq:variational_coupling2}
\end{equation}
Let ${\{\phi^{\ell-1}_j\}_{j \in \pN_{\ell-1}}}$ be a basis of $\pV_{\ell-1}$, ${\{\phi^\ell_k\}_{k \in \pN_{\ell}}}$ be a basis of $\pV_\ell$ and ${\{\theta^\ell_i\}_{i \in \pN_\mu}}$ be basis of the multiplier space $\pM_\ell$, where $\pN_{\ell-1}$, $\pN_{\ell}$ and $\pN_\mu$ denote the set of nodes associated with respective FE space.

Writing the functions $v \in \pV_{\ell-1}$ and $w \in \pV_{\ell}$ as a linear combination of the basis functions, we get $v = \sum_{j \in \pN_{\ell-1}} v_j \phi^{\ell-1}_j$ and $w = \sum_{k \in \pN_{\ell}} w_k \phi^{\ell}_k$, with coefficients ${\{v_j\}_{j \in \pN_{\ell-1}}}$ and ${\{w_k\}_{k \in \pN_{\ell}}}$.
Now inserting the respective basis function in \eqref{eq:variational_coupling2}, we get
\begin{equation}
  \sum_{j \in \pN_{\ell-1}} v_j (\theta_i^\ell, \phi_j^{\ell-1})_{L^2(I)}  = \sum_ {k \in \pN_{\ell}} w_k (\theta_i^\ell, \phi_k^\ell)_{L^2(I)}, \quad \forall i \in \pN_\mu.
  \label{eq:reformulation}
\end{equation}
If we write the above formulation in matrix notation, a $\bs{B}$ matrix is defined between a fine and a coarse space, with entries $B_{ij}=(\theta_i^\ell, \phi_j^{\ell-1})_{L^2(I)}$ and $\bs{D}$ matrix is defined on a fine level, with entries $D_{ik}=(\theta_i^\ell, \phi_k^\ell)_{L^2(I)}$.
The formulation \eqref{eq:reformulation} in matrix-vector form is given as
\begin{equation*}
  \bs{B}\bv = \bs{D} \bs{w}.
  \label{eq:disctrete_form}
\end{equation*}
As Lagrange multiplier space and the finite element space on fine level $\ell$ are of the same dimension, $\bs{D} \in \bR^{\vert \pN_\ell \vert \times \vert \pN_\mu \vert}$ is a square matrix, while $\bs{B} \in \bR^{\vert \pN_\ell \vert \times \vert  \pN_{\ell-1} \vert }$ is a rectangular matrix, here \(\vert \cdot \vert\) denotes the cardinality of a given set.
The vectors $\bv$ and $\bs{w}$ are representation of function $v,w$ on level $\ell-1$ and $\ell$, respectively.
The formula for computing the transfer operator can be expressed algebraically as
\begin{equation}
  \bs{w}=\bs{D}^{-1} \bs{B}\bv = \bs{T}\bv.
  \label{eq:projection}
\end{equation}
The matrix ${\bs{T} \in \bR^{\vert \pN_\ell \vert \times \vert \pN_{\ell-1} \vert }}$ is the discrete representation of the projection operator $\Pi_{\ell-1}^{\ell}$, which now we will use as a prolongation operator in the multigrid method.

\subsubsection{$L^2$-projections}
The choice of different multiplier spaces in the above formulation can lead to different transfer operators.
The Lagrange multiplier space can be the same as the finite element space $ \pM_\ell = \pV_\ell$.
Thus, we take the same basis functions $\{\phi^\ell_i\}_{i \in \pN_\mu}$ for $\pM_\ell$ and $\pV_\ell$.
In this particular case, the scaled mass matrix $\bs{B}$, between a coarse and a fine space, has the elements ${B_{ij}= (\phi_i^\ell, \phi_j^{\ell-1})_{L^2(I)}}$.
The $\bs{D}$ matrix is the mass matrix on a fine level with the elements ${D_{ik} = (\phi_i^\ell ,\phi_k^\ell)_{L^2(I)}}$.

The usage of the transfer operator computed with the $L^2$-projection does not guarantee a computationally efficient multigrid algorithm.
In the multigrid algorithm, we employ the Galerkin assembly approach \eqref{eq:galerkin_assembly} to create coarse level operators, hence it is necessary to compute the transfer operator $\bs{T}$.
For the computation of the transfer operator, a sparse block diagonal matrix $\bs{D}$ has to be inverted.
As the inverse of matrix $\bs{D}$ is dense, it also results in a denser transfer operator which gives rise to coarse spaces with the basis functions with global support.
Hence, the computation of the coarse level quantities and the Galerkin assembly becomes more expensive as the matrix-vector multiplications have to be performed on the dense system.

\subsubsection{Pseudo-$L^2$-projections}
In order to make the transfer operator sparse and to reduce the computational complexity of the application of the transfer operator, we choose a different multiplier space in \eqref{eq:variational_coupling}~\cite{wohlmuth_mortar_2000,dickopf_evaluating_2014}.
The basis functions which span the new multiplier space are chosen in such a way that they are biorthogonal to the standard Lagrange FE basis with respect to $L^2$-inner product.

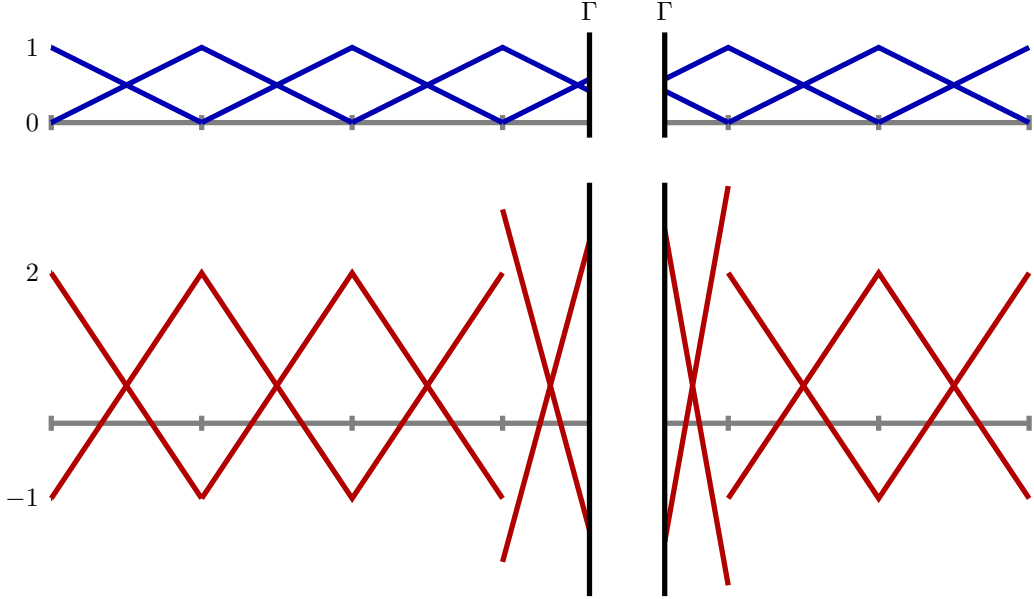
\begin{figure*}[t]
  \centering
  \begin{tikzpicture}[scale=1]
    \begin{scope}[shift={(-0.5,1)}]
  \foreach \i in {4,6,...,11}{
    \draw[gray,line width=0.7mm](\i,-0.1) -- (\i,0.1);
  }

  \draw[line width=0.7mm, gray](4,0) -- (11.1556,0);

  \draw[line width=0.7mm] (4,1)--(4,1) node [left = 0.01 cm] {$1$};
  \draw[line width=0.7mm] (4,0)--(4,0) node [left = 0.01 cm] {$0$};
  \draw[blue!70!black,line width=0.7mm](4,1) -- (6,0);
  \draw[blue!70!black,line width=0.7mm](4,0) -- (6,1) -- (8,0);
  \draw[blue!70!black,line width=0.7mm](6,0) -- (8,1) -- (10,0);
  \draw[blue!70!black,line width=0.7mm](8,0) -- (10,1) -- (11.1556,0.4222);
  \draw[blue!70!black,line width=0.7mm](10,0) -- (11.1556,0.5778);

  \draw[line width=0.7mm] (11.1556,-0.2)--(11.1556,1.2) node [above = 0.01 cm] {$\Gamma$} ;
\end{scope}

\begin{scope}[shift={(0.5,1)}]
  \foreach \i in {12,14,...,16}{
    \draw[gray,line width=0.7mm](\i,-0.1) -- (\i,0.1);
  }
  \draw[line width=0.7mm,gray](11.1556,0) -- (16,0);

  \draw[blue!70!black,line width=0.7mm](11.1556,0.4222) -- (12,0);
  \draw[blue!70!black,line width=0.7mm](11.1556,0.5778) -- (12,1) -- (14,0);
  \draw[blue!70!black,line width=0.7mm](12,0) -- (14,1) -- (16,0);
  \draw[blue!70!black,line width=0.7mm](14,0) -- (16,1);

  \draw[line width=0.7mm] (11.1556,-0.2)--(11.1556,1.2) node [above = 0.01 cm] {$\Gamma$} ;
\end{scope}

\begin{scope}[shift={(-0.5,-3)}]
  \foreach \i in {4,6,...,11}{
    \draw[gray,line width=0.7mm](\i,-0.1) -- (\i,0.1);
  }

  \draw[line width=0.7mm] (4,2)--(4,2) node [left = 0.01 cm] {$2$};
  \draw[line width=0.7mm] (4,-1)--(4,-1) node [left = 0.01 cm] {$-1$};

  \draw[gray, line width=0.7mm](4,0) -- (11.1556,0);

  \draw[line width=0.7mm,red!70!black](4,2) -- (6,-1);
  \draw[line width=0.7mm,red!70!black](4,-1) -- (6,2) -- (8,-1) -- (10,2);
  \draw[line width=0.7mm,red!70!black](6,-1) -- (8,2) -- (10,-1);

  \draw[line width=0.7mm,red!70!black](10,-1.8444) -- (11.1556,2.4222);
  \draw[line width=0.7mm,red!70!black](10,2.8444) -- (11.1556,-1.4222);

  \draw[line width=0.7mm] (11.1556,3.2)--(11.1556,-2.3);
\end{scope}

\begin{scope}[shift={(0.5,-3)}]
  \foreach \i in {12,14,...,16}{
    \draw[gray,line width=0.7mm](\i,-0.1) -- (\i,0.1);
  }
  \draw[gray,line width=0.7mm](11.1556,0) -- (16,0);
  \draw[red!70!black,line width=0.7mm](11.1556,2.5778) -- (12,-2.1556);
  \draw[red!70!black,line width=0.7mm](11.1556,-1.5778) -- (12,3.1556);

  \draw[red!70!black,line width=0.7mm](12,-1) -- (14,2) -- (16,-1);
  \draw[red!70!black,line width=0.7mm](12,2)--(14,-1) -- (16,2);

  \draw[line width=0.7mm] (11.1556,3.2)--(11.1556,-2.3);
\end{scope}

   \end{tikzpicture}
  \caption[biorthogonal basis]{Lagrange basis functions (in blue) for the cut domain and corresponding biorthogonal basis functions (in red), the biorthogonal basis for the cut-element have been modified to satisfy biorthogonality condition for truncated support}
  \label{fig:biorth_basis_l2}
\end{figure*}

We define the dual space, $\pM_\ell:=\text{span}\{\psi_i^\ell\}_{i \in \pN_\mu}$, where $\psi_i^\ell$ are defined as the dual functions which satisfy the following biorthogonality condition on domain $\omega$
\[
  (\psi^\ell_i,\phi^\ell_j)_{L^2(\omega)}= \delta_{ij} (\phi_i^\ell,\mathbbm{1})_{L^2(\omega)}, \quad \forall i,j \in \pN_\mu = \pN_{\ell}.
  \label{eq:dual_product}
\]
For linear and bilinear elements, it is possible to compute the biorthogonal basis ($\psi_j$) as a linear combination of the Lagrange basis ($\phi_k$) for each element $K \in \widetilde{\pT}_\ell$,
\[
  \psi_j = C_{jk} \phi_k,
\]
where $C_{jk}$ denotes the entries of the matrix $\bs{C}$, which represents the coefficients of the linear combination.
Using this information, we can compute the coefficients matrix $\bs{C}_K$, for each element $K \in \pT_{\ell,i}$
\[
  C_{jk} (\phi^\ell_j, \phi^\ell_k)_{L^2(K_i)} =  \delta_{jk} (\phi_j^\ell,\mathbbm{1})_{L^2(K_i)}, \quad \forall j,k \in \pN_K,
\]
where $\pN_K$ denotes set of nodes of a given triangulation $K$.
The matrix representation of the multiplication of element-wise basis function is defined as $\widetilde{M}_{jk}^K =  (\phi^\ell_j, \phi^\ell_k)_{L^2(K_i)}$ and ${\widetilde{D}_{jk}^K = \delta_{jk} (\phi_j^\ell,\mathbbm{1})_{L^2(K_i)}}$.
Here, the elemental matrix $\widetilde{\bs{D}}_K$ is diagonal and $\widetilde{\bs{M}}_K$ is an elemental mass matrix with the entries $\widetilde{D}_{jk}^K$ and  $\widetilde{M}^K_{jk}$, respectively.
In the algebraic form, it can be written as
\[
  \bs{C}_K \widetilde{\bs{M}}_K = \widetilde{\bs{D}}_K.
\]
The element-wise coefficients of the linear combination could be computed as
\begin{equation}
  \bs{C}_K = \widetilde{\bs{D}}_K (\widetilde{\bs{M}}_K)^{-1}.
\end{equation}
As the Lagrange basis functions of the cut elements have truncated support, the biorthogonal basis functions computed using this strategy are not necessarily continuous or regular for the cut elements.
In Figure~\ref{fig:biorth_basis_l2}, we can see the Lagrange basis and the corresponding biorthogonal basis functions for a cut mesh.

In the discrete setting, the matrix entries $D_{ik} = (\psi_i^\ell, \phi^\ell_k)_{L^2(I)}$ computed using this strategy is diagonal.
The diagonal mass matrix $\bs{D}$ is computationally trivial to invert, and the inverse of the matrix is also diagonal.
The matrix $\bs{B}$, defined between a coarse and a fine level can be given with the entries ${B_{ij} = (\psi_i^\ell, \phi^{\ell-1}_j)_{L^2(I)}}$.
The transfer operator computed using this method has a sparse structure, and the support of the basis functions on the coarse level is also localized compared to the standard $L^2$-projection operator.

\begin{remark}
  Computation of the coefficient matrix $\bs{C}_K$ is only necessary for the cut elements as the support of the basis function is arbitrary. For the uncut elements, we compute the coefficient matrix for any element and reuse it for all others.
\end{remark}

\subsection{Semi-geometric Multigrid Method}
We employ the pseudo-$L^2$-projection based transfer operator in our semi-geometric multigrid method.
By construction, it is clear that the transfer operator for the XFEM discretization treats both subdomains separately.
It can be regarded as an additive subspace splitting strategy to compute the transfer operator.
This ensures that the information transfer between levels is restricted to each subdomain, and there is no cross-information transfer across the interface.
This is quite essential for ensuring the robustness of the multigrid method.

Let $\{\lambda_i^{\ell}\}_{i\in\pN_{\ell}}$ be the basis of a FE space $\pX_\ell$ and let $\{\lambda_j^{\ell-1}\}_{j\in\pN_{\ell-1}}$ be the basis of a coarse level FE space, $\pX_{\ell-1}$.
The pseudo-$L^2$-projection based prolongation operator can be used to compute the basis function on the coarse levels,
\[
  \lambda_i^{\ell-1}:=\sum_{i\in\pN_{\ell}} (\bs{T}_{\ell-1}^{\ell})_{ij}\lambda_{j}^{\ell}, \quad \forall j\in \pN_{\ell-1}.
\]
Thus, the basis functions at the coarse level can be computed as a linear combination of the basis functions from the fine level.
As we can see in the Figure~\ref{fig:multilevel_decomp}, the basis functions at the coarse level are piecewise linear with respect to the finest level mesh, given the basis functions are linear also on the finest level.

Considering the abstract weak formulation of all Nitsche-based XFEM discretizations, we rewrite the problem agnostic of the stabilization method as, find $u_h \in \pV_L=\pX_L$, such that:
\[
  a(u_h,v_h) = F(v_h), \quad \forall v_h \in \pV_L.
\]
We denote the algebraic representation of the bilinear form as $\bs{A}_L$, also known as a stiffness matrix, with the entries $A^L_{ij}:=a(\lambda_j^{L},\lambda_i^L)$ for all $i,j\in \pN_L$.
The right hand side is represented as $\bs{f}_L$, with the local entries $(f_L)_i := F(\lambda_i^L)$, for all $i \in \pN_L$.
We can write the algebraic variant of our problem as
\begin{equation}
  \bs{A}_L \bu = \bs{f}_L,
  \label{eq:disctrete_eq}
\end{equation}
where unknown $\bu$ represents a vector of the coefficients $\bu = (u_i)_{i\in\pN_L}$, which are associated with the finite element approximation on space $\pX_L$, $u =\sum_{i \in \pN_L} u_i \lambda_i^L$.

From the practical aspect, we need a preliminary step before the multigrid algorithm can be invoked.
In the preparation step, we compute a sequence of prolongation operators $(\bs{T}_{\ell-1}^\ell)_{\ell\in \{1,\dots,L\}}$ between successive levels in the hierarchy of FE spaces.
These prolongation operators are used to project a coarse level correction to a fine level and the adjoint of these operators are used as restriction operators to project the residual from a fine level to a coarse level.
The next stage of the setup consists of computation of the coarse level stiffness matrices, $(\bs{A}_{\ell})_{\ell\in\{0, \dots, L-1\}}$.
We use the Galerkin assembly approach to compute the coarse level stiffness matrices, defined as
\begin{equation}
  \bs{A}_{\ell-1} := (\bs{T}_{\ell-1}^{\ell})^T \bs{A}_\ell \bs{T}_{\ell-1}^{\ell}, \quad \forall \ell \in \{1,\dots,L\}.
  \label{eq:galerkin_assembly}
\end{equation}
The Galerkin assembly is essential in the semi-geometric multigrid algorithm.
If the assembly of the stiffness matrix is done on each level the coarse level operator would be computed on the FE spaces $(\pV_\ell)_{\ell\in\{0,\dots,L-1\}}$, while the Galerkin assembly assures the stiffness matrix is recursively constructed in the nested FE spaces $(\pX_\ell)_{\ell\in\{0,\dots,L-1\}}$.

\begin{algorithm*}[t]
  \caption{Setup Semi-geometric Multigrid algorithm}
  \SetKwData{Left}{left}\SetKwData{This}{this}\SetKwData{Up}{up}
  \SetKwComment{Comment}{$\triangleright$\ }{}
  \SetKwFunction{Union}{Union}\SetKwFunction{FindCompress}{FindCompress}
  \SetKwInOut{Input}{Input}\SetKwInOut{Output}{Output}
  \Input{$\bs{A}_L$, $(\pT_\ell)_{\ell=0,\dots,L}$}
  \Output{$(\bs{T}_{\ell-1}^{\ell})_{\ell = 1,\dots,L},\ (\bs{A}_\ell)_{\ell = 0,\dots,L-1}$}
  Function: Setup SMG\\
  \For {$\ell \gets L,\dots,1$}
  {
    ${\bs{T}_{\ell-1}^{\ell}} \mapsfrom (\pT_{\ell},\pT_{\ell-1}) $ \Comment*[r]{assemble prolongation operator}
    $\bs{A}_{\ell-1} \mapsfrom ({\bs{T}_{\ell-1}^{\ell}})^T \bs{A}_{\ell} \bs{T}_{\ell-1}^{\ell}$ \Comment*[r]{coarse level assembly}
  }
\end{algorithm*}
\begin{algorithm*}[t]
  \caption{Semi-geometric Multigrid algorithm - $V$-cycle}
  \SetKwData{Left}{left}\SetKwData{This}{this}\SetKwData{Up}{up}
  \SetKwComment{Comment}{$\triangleright$\ }{}
  \SetKwFunction{Union}{Union}\SetKwFunction{FindCompress}{FindCompress}
  \SetKwInOut{Input}{Input}\SetKwInOut{Output}{Output}
  \Input{$(\bs{A}_\ell)_{\ell=0,\dots,L}$, $\bs{r}_L$, $L$, $\nu_1$, $\nu_2$, $(\bs{T}_{\ell-1}^\ell)_{\ell=1,\dots,L}$}
  \Output{$\bs{c}_L$}
  \BlankLine
  Function: SMG($\bs{A}_{\ell}$, $\bs{r}_\ell$, $\ell$, $\nu_1$, $\nu_2$, $\bs{T}_{\ell-1}^\ell$):\\
  \eIf{$\ell \neq 0$}{
  $\bs{c}_{\ell}   \mapsfrom \bs{0}$                                                     \Comment*[r]{initialize correction}
  $\bs{c}_{\ell}   \mapsfrom $ Smoother($\bs{A}_\ell$, $\bs{c}_\ell$, $\bs{r}_\ell$, $\nu_1$)    \Comment*[r]{$\nu_1$ pre-smoothing steps}
  $\bs{r}_{\ell-1} \mapsfrom (\bs{T}_{\ell-1}^{\ell})^T (\bs{r}_\ell - \bs{A}_\ell \bs{c}_\ell)$  \Comment*[r]{restriction}
  $\bs{c}_{\ell-1} \mapsfrom $ SMG($\bs{A}_{\ell-1}$, $\bs{r}_{\ell-1}$, $\ell-1$, $\nu_1$, $\nu_2$, $\bs{T}_{\ell-2}^{\ell-1}$) \Comment*[r]{coarse level cycle}
  $\bs{c}_{\ell}   \mapsfrom \bs{c}_\ell$ + $\bs{T}_{\ell-1}^{\ell}$ $\bs{c}_{\ell-1}$                     \Comment*[r]{prolongation}
  $\bs{c}_{\ell}   \mapsfrom $ Smoother($\bs{A}_\ell$, $\bs{c}_\ell$, $\bs{r}_\ell$, $\nu_2$)    \Comment*[r]{$\nu_2$ post-smoothing steps}
  }
  {
  $\bs{c}_0 \mapsfrom  \bs{A}_0^{-1} \bs{r}_0 $                             \Comment*[r]{direct solver}
  }
\end{algorithm*}

After the computation of the coarse level stiffness matrices and the prolongation operators, we invoke the multigrid algorithm.
The multigrid iterations can be used for preconditioning or solving a linear system, so the SMG can be written in an abstract way to return the correction $\bs{c}_L$, rather than the iterate explicitly.
The residual on the finest level is given as
\[
  \bs{r}_L = \bs{f}_L - \bs{A}_L \bu.
\]
The $V$-cycle semi-geometric multigrid algorithm is given in Algorithm 2, where $\nu_1$, $\nu_2$ are number of pre-smoothing and post-smoothing steps respectively.

 \section{Numerical Results} \label{sect:results}
In this section, we introduce different examples for the numerical experiments.
We compare the different variants of Nitsche's methods on these examples and compare the discretization errors and the condition number of the system matrices.
Also, we evaluate the performance of the semi-geometric multigrid method as a solution method and as a preconditioner for the same examples.
In this section for brevity, we drop the subscript $L$ from the stiffness matrix ${\bA=\bA_L}$.

\subsection{Problem Description}
We consider a domain $\Omega = [0,1]^2$ with two different types of interfaces, a linear interface, and a circular interface.
All the experiments are carried out on a triangular structured mesh as shown in Figure~\ref{fig:multilevel_decomp}.
We start with a mesh that has $100$ elements in each direction, denoted as $L1$ and uniformly refine the mesh $L1$ to obtain different meshes as shown in Table~\ref{tab:mesh_hierarchy}.
We use the same mesh hierarchy to measure the discretization error and later as the multilevel hierarchy in the multigrid method.
We remark that the dofs and number of elements already include the enriched nodes and the enriched elements.
\begin{table*}[t]\centering
  \begin{tabular}{ c | c |r  r | r  r  }
    \multirow{2}{*}{levels}                          &
    \multirow{2}{*}{$h_{\max}$}                      &
    \multicolumn{2}{c|}{Linear Interface $\Gamma_l$} &
    \multicolumn{2}{c}{Circular Interface $\Gamma_c$}                                                                          \\
    \cline{3-6}
                                                     &             & dofs      & No. of Elements & dofs      & No. of Elements \\\hline\hline
    $L1$                                             & 1.41421E-02 & 10,403    & 20,200          & 10,767    & 20,566          \\\hline
    $L2$                                             & 7.07107E-03 & 40,803    & 80,400          & 41,531    & 81,130          \\\hline
    $L3$                                             & 3.53553E-03 & 161,603   & 320,800         & 163,063   & 322,262         \\\hline
    $L4$                                             & 1.76777E-03 & 643,203   & 1,281,600       & 646,127   & 1,284,526       \\\hline
    $L5$                                             & 8.83883E-04 & 2,566,403 & 5,123,200       & 2,572,251 & 5,129,050       \\
  \end{tabular}
  \caption {The multilevel hierarchy}
  \label{tab:mesh_hierarchy}
\end{table*}
We consider problems with continuous and discontinuous coefficients to analyze the effect of different variants of Nitsche's method on the condition number of the linear system and the numerical accuracy of the discretization.

\paragraph{Example 1}\label{para:example1}
Consider a Poisson problem where $\alpha_1 = \alpha_2 = 1$, and a linear interface $\Gamma_l$.
For this example, the right-hand side $f_1$ and the Dirichlet boundary condition are chosen in such a way that the exact solution, $u_1 = (\exp(-500s)-1)(\exp(-500t)-1)(\exp(-500\hat{y})-1)(1-3\hat{r})^2$ is satisfied.
Here, $s:= (x-1/3)^2$, $t:=(x-2/3)^2$, $\hat{x}:=(x-1/2)^2$, $\hat{y}:=(y-1/2)^2$ and $\hat{r}:=\hat{x}+\hat{y}$.
The linear interface $\Gamma_l$ is defined as a zero level set of the function $\Lambda_l(x,y) := x - 1/\sqrt{2}$.
We have deliberately chosen the location of interface $\Gamma_l$ in such a way that the interface would stay close to edges of elements for all levels of the refinement and would not coincide with the element edges.
Thus, the enriched elements are divided into disproportional fractions.

\paragraph{Example 2}\label{para:example2}
For this example, we consider a problem with discontinuous coefficients and a circular interface $\Gamma_c$.
The circular interface $\Gamma_c$ is defined as a zero level set of a function $\Lambda_c(\bx) := r_0^2-\|\bx-\bs{c}\|_2^2$ with radius, $r_0^2 = {3-\sqrt{2}}$, and $\bs{c}$ is the center of the circle $(0.5,0.5)$.
The circular interface decomposes the domain $\Omega$ into $\Omega_1$, where $\Lambda_c(\bx) > 0$ and $\Omega_2$ where $\Lambda_c(\bx) < 0$.

We consider a Poisson problem, where we choose coefficients as ${\alpha_1 = \{10^{-1},10^{-5},10^{-9}\}}$ and ${\alpha_2 = 1}$.
For this example, the right-hand side is chosen as ${f_2 = -4\alpha_1 \alpha_2}$, and the Dirichlet boundary conditions satisfy the exact solution,
\begin{equation*}
  u_2(\bx) =
  \begin{cases}
    \alpha_2(\|\bx-\bs{c}\|_2^2 - r_0^2),  & \text{ if } \bx \in \Omega_1, \\
    \alpha_1(\|\bx-\bs{c}\|_2^2 - r_0^2) , & \text{ if } \bx \in \Omega_2.
  \end{cases}
\end{equation*}

\paragraph{Example 3}\label{para:example3}
Here we consider the same circular interface as in \nameref{para:example2}.
In this example, we consider the coefficients as ${\alpha_1=1}$ and ${\alpha_2=\{10,10^5,10^9\}}$.
The right-hand side is chosen as, $f_3 = -4$, and the Dirichlet boundary conditions satisfy the exact solution,
\begin{equation*}
  u_3(\bx) =
  \begin{cases}
    \dfrac{\|\bx-\bs{c}\|_2^2}{\alpha_1},                                   & \text{ if } \bx \in \Omega_1, \\
    \dfrac{\|\bx-\bs{c}\|_2^2 - r_0^2}{\alpha_2} + \dfrac{r^2_0}{\alpha_1}, & \text{ if } \bx \in \Omega_2.
  \end{cases}
\end{equation*}

\subsection{Comparison of Variants of Nitsche's Methods}
We compare the convergence rate of error in the $L^2$-norm and the mesh-dependent energy norm against the condition number of the stiffness matrix, denoted as $\kappa(\bA)$.
The motivation behind this comparison is to investigate the influence of previously discussed variants of Nitsche's methods on condition numbers and the discretization errors.
As the background mesh and the location of the interface is fixed for a given example, the condition number of the system matrix is only affected by the choice of stabilization parameter $\gamma$ and the weighting parameters $\beta_i$.

\subsubsection{Discretization Error}
From Figure~\ref{fig:convergence_1} and Figure \ref{fig:convergence_23}, it is clear that all variants of Nitsche's methods have almost identical rate of convergence for the discretization error in both norms.
As the error estimates in \eqref{eq:xfem_error} suggest, the rate of convergence of discretization error in the $L^2$-norm is ${O}(h^2)$ and in the mesh-dependent energy norm is ${O}(h)$.
Figure~\ref{fig:convergence_1} demonstrates that FEM and XFEM methods have the same approximation properties, since both methods produce the same discretization error for the same mesh size.
The discretization error of \nameref{para:example2} and \nameref{para:example3} is also almost identical in both norms for all variants of Nitsche's method.
Thus, we can conclude that all variants of Nitsche's formulations are robust with respect to continuous and with respect to highly varying coefficients.
\tikzset{external/export next=false}
\begin{figure*}[t]
  \centering
  \begin{tikzpicture}[scale=1]
    \begin{loglogaxis}[ylabel={Error}, xlabel=$\kappa(\bs{A})$, minor tick, grid=major,legend style={font=\small},legend style={at={(axis cs:1e7,1)},anchor=north west},label style={font=\small}, tick label style={font=\small}]
\addplot [very thick,cyan!70!black,mark=diamond*]table[x=FEM_cond, y=FEM_l2, col sep=comma]{example1.csv};
      \label{pgfplots:1FEML2}
      \addplot [very thick,red!70!black,mark=square*]table[x=NXFEMEV_cond, y=NXFEMEV_l2, col sep=comma]{example1.csv};
      \label{pgfplots:1N1L2}
      \addplot [very thick,green!70!black,mark=triangle*]table[x=NXFEMLO_cond, y=NXFEMLO_l2, col sep=comma]{example1.csv};
      \label{pgfplots:1N2L2}
      \addplot [very thick,blue!70!black,mark=*]table[x=NXFEMGP_cond, y=NXFEMGP_l2, col sep=comma]{example1.csv};
      \label{pgfplots:1N3L2}

      \addplot [densely dashed,very thick,cyan!70!black,mark=diamond*,every mark/.append style={solid}]table[x=FEM_cond, y=FEM_sh1, col sep=comma]{example1.csv};
      \label{pgfplots:1FEMH1}
      \addplot [densely dashed,very thick,red!70!black,mark=square*,every mark/.append style={solid}]table[x=NXFEMEV_cond, y=NXFEMEV_brh1, col sep=comma]{example1.csv};
      \label{pgfplots:1N1BH1}
      \addplot [densely dashed,very thick,green!70!black,mark=triangle*,every mark/.append style={solid}]table[x=NXFEMLO_cond, y=NXFEMLO_brh1, col sep=comma]{example1.csv};
      \label{pgfplots:1N2BH1}
      \addplot [densely dashed,very thick,blue!70!black,mark=*,every mark/.append style={solid}]table[x=NXFEMGP_cond, y=NXFEMGP_brh1, col sep=comma]{example1.csv};
      \label{pgfplots:1N3BH1}
    \end{loglogaxis}
    \matrix [ draw, matrix of nodes, anchor = west, node font=\small,
    column 1/.style={nodes={align=left,text width=1cm}},
    column 2/.style={nodes={align=center,text width=2cm}},
    column 3/.style={nodes={align=center,text width=2cm}}
    ] at (current axis.east)
    {
    & $\|u - u_h\|_{L^2(\Omega)}$  & $\vertiii{u-u_h}_h$    \\
    FEM   & \ref{pgfplots:1FEML2}  &  \ref{pgfplots:1FEMH1} \\
    \ref{eq:NXFEM-EV}  & \ref{pgfplots:1N1L2}   &  \ref{pgfplots:1N1BH1} \\
    \ref{eq:NXFEM-LO}  & \ref{pgfplots:1N2L2}   &  \ref{pgfplots:1N2BH1} \\
    \ref{eq:NXFEM-GP}  & \ref{pgfplots:1N3L2}   &  \ref{pgfplots:1N3BH1} \\
    };

  \end{tikzpicture}
  \caption{Convergence of error in the $L^2$-norm and the mesh-dependent energy norm (except for the FEM formulation, where we use the $H^1$-seminorm) for the different variants of Nitsche's method applied to \nameref{para:example1}}
  \label{fig:convergence_1}
\end{figure*}
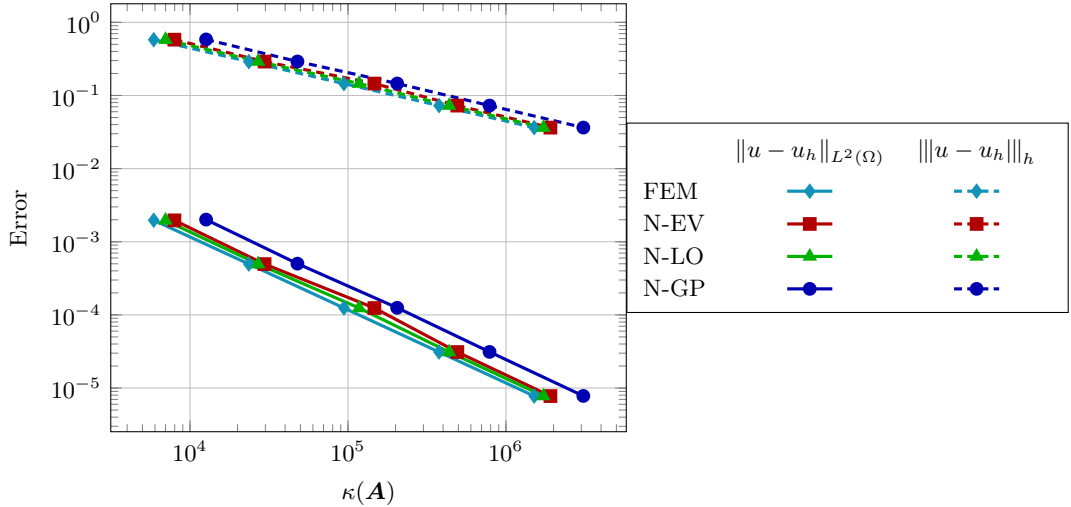
\tikzset{external/export next=false}
\begin{figure*}[t]
  \centering
  \begin{tikzpicture}[]
    \begin{groupplot}[
        group style={
group size = 3 by 1,
x descriptions at=edge bottom,
horizontal sep=1pt,
          },
scale=0.78,ymode=log,xmode=log,xlabel=$\kappa(\bs{A})$, minor tick num=4, grid=major,
label style={font=\small}, tick label style={font=\small}, legend style={font=\small}
      ]
\nextgroupplot[align=left, title={{\small \nameref{para:example2}: $\alpha_1=10^{-1},\alpha_2=1$ \\ \small \nameref{para:example3}: $\alpha_1=1, \quad \ \ \alpha_2=10$} },ylabel={Error}]
      \addplot [very thick,red!70!black,mark=square*]table[x=NXFEMEV_cond, y=NXFEMEV_l2, col sep=comma]{example2alpha1.csv};
      \label{pgfplots:2N1L2}
      \addplot [very thick,green!70!black,mark=triangle*]table[x=NXFEMLO_cond, y=NXFEMLO_l2, col sep=comma]{example2alpha1.csv};
      \label{pgfplots:2N2L2}
      \addplot [very thick,blue!70!black,mark=*]table[x=NXFEMGP_cond, y=NXFEMGP_l2, col sep=comma]{example2alpha1.csv};
      \label{pgfplots:2N3L2}

      \addplot [densely dashed,very thick,red!70!black,mark=square*,every mark/.append style={solid}]table[x=NXFEMEV_cond, y=NXFEMEV_brh1, col sep=comma]{example2alpha1.csv};
      \label{pgfplots:2N1BH1}
      \addplot [densely dashed,very thick,green!70!black,mark=triangle*,every mark/.append style={solid}]table[x=NXFEMLO_cond, y=NXFEMLO_brh1, col sep=comma]{example2alpha1.csv};
      \label{pgfplots:2N2BH1}
      \addplot [densely dashed,very thick,blue!70!black,mark=*,every mark/.append style={solid}]table[x=NXFEMGP_cond, y=NXFEMGP_brh1, col sep=comma]{example2alpha1.csv};
      \label{pgfplots:2N3BH1}

\nextgroupplot[align=left, title={{\small \nameref{para:example2}: $\alpha_1=10^{-5},\alpha_2=1$ \\ \small \nameref{para:example3}: $\alpha_1=1, \quad \ \ \alpha_2=10^5$}},yticklabels={}]
      \addplot [very thick,red!70!black,mark=square*]table[x=NXFEMEV_cond, y=NXFEMEV_l2, col sep=comma]{example2alpha5.csv};
      \addplot [very thick,green!70!black,mark=triangle*]table[x=NXFEMLO_cond, y=NXFEMLO_l2, col sep=comma]{example2alpha5.csv};
      \addplot [very thick,blue!70!black,mark=*]table[x=NXFEMGP_cond, y=NXFEMGP_l2, col sep=comma]{example2alpha5.csv};

      \addplot [densely dashed,very thick,red!70!black,mark=square*,every mark/.append style={solid}]table[x=NXFEMEV_cond, y=NXFEMEV_brh1, col sep=comma]{example2alpha5.csv};
      \addplot [densely dashed,very thick,green!70!black,mark=triangle*,every mark/.append style={solid}]table[x=NXFEMLO_cond, y=NXFEMLO_brh1, col sep=comma]{example2alpha5.csv};
      \addplot [densely dashed,very thick,blue!70!black,mark=*,every mark/.append style={solid}]table[x=NXFEMGP_cond, y=NXFEMGP_brh1, col sep=comma]{example2alpha5.csv};

\nextgroupplot[align=left, title={{\small \nameref{para:example2}: $\alpha_1=10^{-9},\alpha_2=1$ \\ \small \nameref{para:example3}: $\alpha_1=1, \quad \ \ \alpha_2=10^9$}},yticklabels={}]
      \addplot [very thick,red!70!black,mark=square*]table[x=NXFEMEV_cond, y=NXFEMEV_l2, col sep=comma]{example2alpha9.csv};
      \addplot [very thick,green!70!black,mark=triangle*]table[x=NXFEMLO_cond, y=NXFEMLO_l2, col sep=comma]{example2alpha9.csv};
      \addplot [very thick,blue!70!black,mark=*]table[x=NXFEMGP_cond, y=NXFEMGP_l2, col sep=comma]{example2alpha9.csv};

      \addplot [densely dashed,very thick,red!70!black,mark=square*,every mark/.append style={solid}]table[x=NXFEMEV_cond, y=NXFEMEV_brh1, col sep=comma]{example2alpha9.csv};
      \addplot [densely dashed,very thick,green!70!black,mark=triangle*,every mark/.append style={solid}]table[x=NXFEMLO_cond, y=NXFEMLO_brh1, col sep=comma]{example2alpha9.csv};
      \addplot [densely dashed,very thick,blue!70!black,mark=*,every mark/.append style={solid}]table[x=NXFEMGP_cond, y=NXFEMGP_brh1, col sep=comma]{example2alpha9.csv};

    \end{groupplot}
    \matrix [ draw, matrix of nodes, anchor = north, node font=\small,
    column 1/.style={nodes={align=left,text width=2cm}},
    column 2/.style={nodes={align=center,text width=1.3cm}},
    column 3/.style={nodes={align=center,text width=1.3cm}},
    column 4/.style={nodes={align=center,text width=1.3cm}},
    ] at ($(group c2r1) + (0,-3.2)$)
    {
    &\ref{eq:NXFEM-EV}      & \ref{eq:NXFEM-LO}     & \ref{eq:NXFEM-GP}              \\
    $\|u - u_h\|_{L^2(\Omega)}$ & \ref{pgfplots:2N1L2}  & \ref{pgfplots:2N2L2}  & \ref{pgfplots:2N3L2}  \\
    $\vertiii{u-u_h}_h$       & \ref{pgfplots:2N1BH1} & \ref{pgfplots:2N2BH1} & \ref{pgfplots:2N3BH1} \\
    };

  \end{tikzpicture}
  \caption{Convergence of error in the $L^2$-norm and the mesh-dependent energy norm for the different variants of Nitsche's method applied to \nameref{para:example2} and \nameref{para:example3}}
  \label{fig:convergence_23}
\end{figure*}
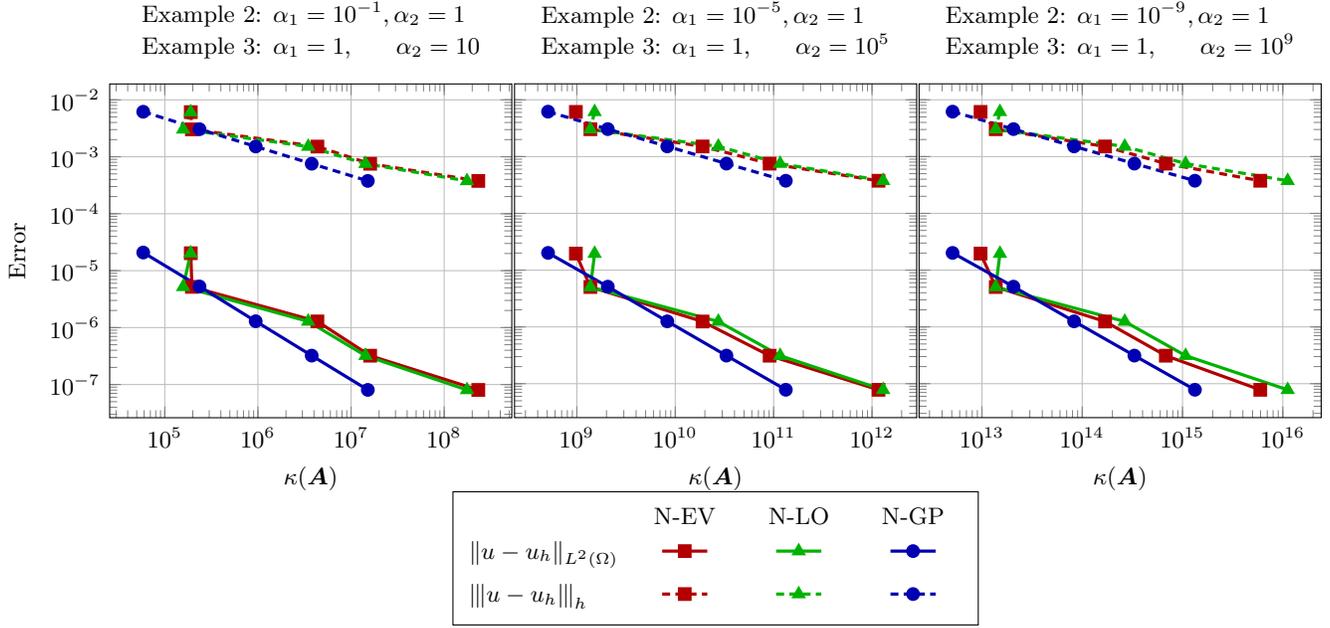

\subsubsection{Condition Numbers}
As the \nameref{para:example1} has continuous coefficients, we can compare the FEM discretization with the XFEM discretization on the same background mesh.
We observe that FEM discretization has smallest $\kappa(\bA)$ for all different mesh sizes in comparison with its XFEM counterparts.
The XFEM discretizations produce the system of linear equations with larger condition numbers in all cases.
The condition number for Nitsche's formulation equipped with the lifting operator is closest to the FEM discretization, while \eqref{eq:NXFEM-EV} formulation is the close second.
The condition numbers of the \eqref{eq:NXFEM-GP} formulation is the largest for all mesh sizes, but this can be attributed to the chosen value of the stabilization parameter.

\nameref{para:example2} and \nameref{para:example3} are different in terms of the coefficients.
The ratio between the smallest and the largest coefficient is kept the same in both examples.
For a given variant of Nitsche's formulation and a given ratio between coefficients, we witness the identical results in terms of error and condition number of the system for both examples.
Application of the Dirichlet boundary condition to the stiffness matrix in both examples causes the distribution of the eigenvalues in the spectrum to vary, but the ratio between the largest and the smallest eigenvalues stays the same.
Figure~\ref{fig:convergence_23} shows a comparison of the condition numbers against the discretization errors for both examples for various coefficients.
It is evident from Figure~\ref{fig:convergence_23}, that the condition numbers and the discretization errors in the $L^2$-norm and the mesh-dependent energy norm for both examples are identical.
The condition number is the smallest for the ghost penalty discretization, regardless of the scale of coefficients, and $\kappa(\bA)$ grows with decreasing mesh size, ${O}(h^{-2})$.
The theoretical estimates of Nitsche's method with the ghost penalty term, suggest that $\kappa(\bA)$ is completely independent of the location of the interface on mesh and only depends on the coefficients.
Experimental results also support the theoretical estimates of \eqref{eq:NXFEM-GP} discretization for all examples.

To our knowledge, there are no theoretical bounds established on $\kappa(\bA)$ for \eqref{eq:NXFEM-LO} and \eqref{eq:NXFEM-EV} discretizations.
For different mesh sizes, $\kappa(\bA)$ is larger than the ghost penalty formulation in almost all cases but it grows as the ratio between the coefficients increases.
We also observe the effect of irregular intersection between the interface and the meshes at different levels.

From this discussion, it is clear that Nitsche's method with the ghost penalty term is most stable among all methods.
The other two variants of Nitsche's methods produce the system matrices with larger condition numbers, but still, these variants are stable as the condition numbers do not grow sporadically.

\subsection{Convergence Studies of the Multigrid Method}
In this section, we evaluate the performance of our semi-geometric multigrid method for different variants of Nitsche's method for the discussed examples.
We employ the multigrid methods as a solution method and as a preconditioner and compare their performance against other preconditioners.

Before discussing the results, we need to define a few metrics to compare the different solution methods.
We define scalar energy product $(\cdot,\cdot)_A$ as,
\[
  (\bu,\bv)_A := \bu^T \bA \bv, \quad \forall \bu,\bv \in  \bR^n,
\]
and the induced energy norm is defined as ${\| \cdot \|_A^2:=(\cdot,\cdot)_A}$.
Our examples have highly varying coefficients, hence to compare all solution methods on the same scale we choose a relative residual in the energy norm as a termination criterion, given as
\begin{equation}
  \dfrac{\| \bs{f} - \bA \bu_k \|_A} { \|\bs{f}-\bA \bu_0 \|_A}  < 10^{-12},
  \label{eq:termination}
\end{equation}
where, $\bu_k$ denotes the $k^{\text{th}}$ iterate.
Additionally, we define the asymptotic convergence rate of an iterative solver as
\[
  \rho^* := \dfrac{\|\bu_{k+1} - \bu_{k} \|_A }{\|\bu_{k} -\bu_{k-1} \|_A},
\]
where the iterate $\bu_{k+1}$ satisfies the termination criterion \eqref{eq:termination}.

\subsubsection{Comparison of Preconditioners}
The system of linear equations arising from Nitsche's method is symmetric positive definite (SPD).
The most natural choice of an iterative solver for such problems is the conjugate gradient (CG) method.
Although the CG method has the best rate of convergence amongst all Krylov solvers for SPD systems, in practice preconditioned CG method is used to ensure the fast convergence and, in some cases, to ensure the convergence of the solver up to a certain tolerance.
We use the preconditioned CG method as a solver in our numerical experiments.
We compare Jacobi, symmetric Gauss-Seidel (SGS) and semi-geometric multigrid methods as preconditioners.
\begin{table*}[t]
  \centering
  \begin{subtable}[]{1\linewidth}
    \centering
    \begin{tabular}{ c | c | c | c  | c  c} & CG-Jacobi & CG-SGS & CG-SMG & SMG & ($\rho^*$) \\ \hline \hline
      \ref{eq:NXFEM-EV} & 6363      & 2048   & 8      & 9   & (0.113)    \\
      \ref{eq:NXFEM-LO} & 5365      & 2001   & 7      & 8   & (0.087)    \\
      \ref{eq:NXFEM-GP} & 5396      & 2051   & 9      & 12  & (0.181)    \\ \end{tabular}
    \caption {\nameref{para:example1}}
  \end{subtable}
  \newline \vspace*{0.4 cm}
  \begin{subtable}[]{1\linewidth}
    \centering
    \begin{tabular}{ c | c | c | c | c | c c } &                   & CG-Jacobi & CG-SGS & CG-SMG & SMG     & ($\rho^*$) \\ \hline \hline
      \multirow{3}{*}{\shortstack[l]{$\alpha_1=10^{-1}$                                          \\ $\alpha_2 = 1 $ } }&
      \ref{eq:NXFEM-EV} & 5276              & 1779      & 8      & 10     & (0.157)              \\
                        & \ref{eq:NXFEM-LO} & 4710      & 1724   & 7      & 9       & (0.112)    \\
                        & \ref{eq:NXFEM-GP} & 4510      & 1787   & 7      & 8       & (0.101)    \\ \hline
      \multirow{3}{*}{\shortstack[l]{$\alpha_1=10^{-5}$                                          \\ $\alpha_2 = 1 $ } }&
      \ref{eq:NXFEM-EV} & 4036              & 1626      & 8      & 10     & (0.172)              \\
                        & \ref{eq:NXFEM-LO} & 4027      & 1626   & 8      & 10      & (0.195)    \\
                        & \ref{eq:NXFEM-GP} & 4494      & 1761   & 7      & 7       & (0.092)    \\ \hline
      \multirow{3}{*}{\shortstack[l]{$\alpha_1=10^{-9}$                                          \\ $\alpha_2 = 1 $ } }&
      \ref{eq:NXFEM-EV} & 4139              & 1661      & 8      & 10     & (0.171)              \\
                        & \ref{eq:NXFEM-LO} & 4132      & 1655   & 8      & 10      & (0.195)    \\
                        & \ref{eq:NXFEM-GP} & 4969      & 1675   & 7      & 7       & (0.092)    \\ \end{tabular}
    \caption {\nameref{para:example2}}
  \end{subtable}
  \newline \vspace*{0.4 cm}
  \begin{subtable}[]{1\linewidth}
    \centering
    \begin{tabular}{ c | c | c | c | c | c c } &                   & CG-Jacobi & CG-SGS & CG-SMG & SMG     & ($\rho^*$) \\ \hline \hline
      \multirow{3}{*}{\shortstack[l]{ $\alpha_1=1$                                               \\ $\alpha_2 = 10$}} &
      \ref{eq:NXFEM-EV} & 5192              & 1794      & 7      & 9      & (0.135)              \\
                        & \ref{eq:NXFEM-LO} & 4655      & 1784   & 7      & 7       & (0.031)    \\
                        & \ref{eq:NXFEM-EV} & 4428      & 1803   & 7      & 7       & (0.033)    \\ \hline
      \multirow{3}{*}{\shortstack[l]{ $\alpha_1=1$                                               \\ $\alpha_2 = 10^5$}} &
      \ref{eq:NXFEM-EV} & 3761              & 1635      & 6      & 7      & (0.029)              \\
                        & \ref{eq:NXFEM-LO} & 3752      & 1635   & 6      & 7       & (0.029)    \\
                        & \ref{eq:NXFEM-GP} & 4417      & 1684   & 6      & 7       & (0.029)    \\ \hline
      \multirow{3}{*}{\shortstack[l]{ $\alpha_1=1$                                               \\ $\alpha_2 = 10^9$}} &
      \ref{eq:NXFEM-EV} & 3535              & 1509      & 6      & 7      & (0.029)              \\
                        & \ref{eq:NXFEM-LO} & 3520      & 1509   & 6      & 7       & (0.029)    \\
                        & \ref{eq:NXFEM-GP} & 3941      & 1550   & 6      & 7       & (0.029)    \\ \end{tabular}
    \caption {\nameref{para:example3}}
  \end{subtable}
  \caption{Number of iterations required to reach the predefined tolerance for different preconditioners,
    the last column shows number of required iterations and the asymptotic convergence rate of multigrid method}
  \label{tab:example_precond}
\end{table*}

The experiments are carried out on the system of linear equations with around $2.5\times 10^{6}$ dofs ($L5$).
Our semi-geometric multigrid method is set up with 5-levels, and symmetric Gauss-Seidel is chosen as smoother with 3 pre-smoothing and 3 post-smoothing steps at each level, and we perform a single $V$-cycle as a preconditioner.

Table~\ref{tab:example_precond} shows the number of iterations required by different methods to reach the termination criterion \eqref{eq:termination}.
We observe that the CG method preconditioned with the Jacobi method has the slowest convergence amongst all solvers.
CG method with SGS as a preconditioner is significantly better than the Jacobi preconditioner, the number of iterations is reduced in more than half for most of the problems.
The best performance from all the preconditioners is clearly shown by SMG method.

The rate of convergence of conjugate gradient method depends on the distribution of the spectrum of $\bA$, and the method performs very well if the eigenvalues are clustered in a certain region of the spectrum, rather than uniformly distributed eigenvalues.
Hence, even with the same condition number for the same method, we require a different number of iteration to reach the termination criterion.
Thus, here we show that the CG-SMG method is stable for all discussed discretization methods and coefficients, the number of iterations required for the convergence stay stable.

\subsubsection{Performance as a Solution Method}
By comparing the results of the SMG method, we observe that Nitsche's method with the ghost penalty stabilization term converges fastest for the highly varying coefficients, while it is slowest for the continuous coefficients.
Even though the difference is not significantly high, a few more iterations can be due to the large value of the stabilization parameter.
For the \eqref{eq:NXFEM-EV} and \eqref{eq:NXFEM-LO} the number of iterations to reach the convergence tolerance stays more or less stable.
The multigrid method can be considered quite robust in terms of the asymptotic convergence rates, as for all the experiments we observe $ \rho^* < 0.2 $.
A multigrid method can be interpreted as a Richardson method with SMG as a preconditioner, and the CG method is known to be far superior to the Richardson method.
Hence, we observe that the number of iterations required is smaller in all cases when the semi-geometric multigrid is chosen as a preconditioner than a solution method.

\subsubsection{Level Independence}
In the next part, we evaluate the performance of CG-SMG method for different levels in the multigrid hierarchy.
The finest level is kept the same as in the previous experiments, and the number of levels used in the multilevel hierarchy is changed.
As we use a direct solver on the coarsest level, the coarse level corrections become increasingly accurate as the number of levels is reduced, but the higher number of levels are computationally cheaper as a smaller linear system of equations has to be solved on the coarsest level.
Table~\ref{tab:example3_levels} demonstrates that the number of iterations stay constant regardless of the number of levels used in the multigrid hierarchy.
We see that the change in the ratio between the coefficients does not affect the performance of CG-SMG method.
This result shows the level independence of the multigrid method as a preconditioner.

\begin{table}[t]
  \begin{subtable}[t]{0.5\textwidth}
    \centering
    \begin{tabular}{ c || c || c | c | c | c }
      \multicolumn{2}{c||}{ \# levels} & 2                 & 3 & 4 & 5     \\ \hline \hline
      \multirow{3}{*}{\shortstack[l]{$\alpha_1=1$                          \\ $\alpha_2 = 10$}} &
      \ref{eq:NXFEM-EV}                & 7                 & 7 & 7 & 7     \\
                                       & \ref{eq:NXFEM-LO} & 6 & 6 & 6 & 7 \\
                                       & \ref{eq:NXFEM-GP} & 7 & 7 & 7 & 7 \\ \hline
      \multirow{3}{*}{\shortstack[l]{$\alpha_1=1$                          \\ $\alpha_2 = 10^5$}} &
      \ref{eq:NXFEM-EV}                & 6                 & 6 & 6 & 6     \\
                                       & \ref{eq:NXFEM-LO} & 6 & 6 & 6 & 6 \\
                                       & \ref{eq:NXFEM-GP} & 6 & 6 & 6 & 6 \\ \hline
      \multirow{3}{*}{\shortstack[l]{$\alpha_1=1$                          \\ $\alpha_2 = 10^9$}} &
      \ref{eq:NXFEM-EV}                & 6                 & 6 & 6 & 6     \\
                                       & \ref{eq:NXFEM-LO} & 6 & 6 & 6 & 6 \\
                                       & \ref{eq:NXFEM-GP} & 6 & 6 & 6 & 6
    \end{tabular}
    \caption {Effect of different number of levels in the hierarchy\\for \nameref{para:example3}}
    \label{tab:example3_levels}
  \end{subtable}
  \hfill
  \begin{subtable}[t]{0.45\textwidth}
    \centering
    \begin{tabular}{ c || c | c | c}
      \# interfaces & \ref{eq:NXFEM-EV} & \ref{eq:NXFEM-LO} & \ref{eq:NXFEM-GP} \\ \hline \hline
      1             & 9                 & 8                 & 9                 \\
      2             & 9                 & 8                 & 9                 \\
      4             & 9                 & 8                 & 9                 \\
      6             & 9                 & 8                 & 9                 \\
      8             & 9                 & 8                 & 9                 \\
      10            & 9                 & 8                 & 9
    \end{tabular}
    \caption {Effect of different number of interfaces in the domain}
    \label{tab:example1_interfaces}
  \end{subtable}
  \caption{Number of iterations required to reach the predefined tolerance for conjugate gradient method preconditioned with semi-geometric method with respect to  different levels in the hierarchy and multiple interfaces in the domain}
  \label{tab:temps}
\end{table}

\subsubsection{Multiple Interfaces}
The last set of experiments is done to demonstrate the robustness of the SMG method with respect to the number of interfaces in a domain.
We consider \nameref{para:example1} for this numerical experiment with continuous coefficients.
The finest level is kept the same as in the previous cases, and the multigrid hierarchy consists of 5-levels.
This test is performed for all the discussed variants of Nitsche's methods with multiple interfaces.

The interfaces are represented by zeros level set of the following functions,
\begin{equation*}
  \Lambda_i(x) :=
  \begin{cases}
    x - 0.1\Big(\frac{1}{\sqrt{2}} +i-1 \Big), & \text{for all } i \in \{1,\ldots,5\},  \\
    x + 0.1\Big(\frac{1}{\sqrt{2}}-i\Big),     & \text{for all } i \in \{6,\ldots,10\}.
  \end{cases}
\end{equation*}
All the interfaces are linear, parallel to the original interface $\Gamma_l$.
We start with a single interface and increase up to 10 interfaces in the domain.
From Table~\ref{tab:example1_interfaces}, we can observe that the proposed multigrid method as a preconditioner is stable, as the number of iterations do not change at all with increasing interfaces in the domain.

Thus, we can conclude that our SMG method is a robust solution strategy.
The method is stable for all variants of Nitsche's method with respect to highly varying coefficients, with respect to the number of levels in the multilevel hierarchy and also with respect to the number of interfaces in the domain.
 \section{Conclusion}
In this paper, we reviewed selected strategies to overcome ill-conditioning related to Nitsche's method for the XFEM discretization.
We discussed two different strategies to implicitly estimate the stabilization parameter and the ghost penalty term to improve the robustness of Nitsche's formulation.
Also, we numerically compared the stability of these methods for continuous and highly varying coefficients in terms of discretization error and condition numbers.
We introduced a semi-geometric multigrid method for the unfitted finite element methods and discussed the $L^2$-projection and pseudo-$L^2$-projection approaches to construct the transfer operator for the XFEM discretization.
In the series of experiments, we demonstrated the robustness of our tailored multigrid method with respect to highly varying coefficients and the number of interfaces in a domain.
Additionally, the multigrid method shows level independent convergence rates when applied to variants of Nitsche's methods.

The multigrid method proposed in this work can be used for any unfitted finite element discretization.
In the future, we aim to extend the multigrid method to more complex problems, for example, contact problems and fluid-structure interaction problems in the unfitted FEM framework.
We also aim to implement the $L^2$-projections for the XFEM discretization in the ParMOONoLith library~\cite{moonolithgit,krause_parallel_2016}.
This library can compute the $L^2$-projection on the complex geometries on distributed computing architecture.
In this way, we can extend the multigrid method from this work to the parallel architecture in order to tackle large-scale problems.

\bibliographystyle{spmpsci}      \def\url#1{}

\begin{appendix}
  \section{Coercivity} \label{sect:appendix}
Following the coercivity of the bilinear form \eqref{eq:N1_blinear}, we have
\begingroup
\allowdisplaybreaks
\whencolumns
{
  \begin{align*}
    a(u,u) & = \sum_{i=1}^2 \|\alpha^{\half}\nabla u \|^2_{L^2(\Omega_i)} - 2 (\curlii[]{\alpha \nabla_{\bn}u}, \llbracket u \rrbracket )_{L^2(\Gamma)} + \|\gamma^\half \llbracket u \rrbracket\|^2_{L^2(\Gamma)}                                                                            \\
           & \geqslant  \sum_{i=1}^2 \|\alpha^{\half}\nabla u \|^2_{L^2(\Omega_i)} - \frac{1}{\epsilon_c} \| \curlii[]{\alpha \nabla_{\bn}u} \|_{H^{-\half}(\Gamma),h } -  \epsilon_c \|\llbracket u \rrbracket \|_{H^{\half}(\Gamma),h} + \gamma \|\llbracket u \rrbracket\|^2_{L^2(\Gamma)} \\
           & = \sum_{i=1}^2 \|\alpha^{\half}\nabla u \|^2_{L^2(\Omega_i)} + \Big(\frac{1}{\epsilon_c} - \frac{2}{\epsilon_c} \Big) \| \curlii[]{\alpha \nabla_{\bn}u} \|_{H^{-\half}(\Gamma),h }
    + \sum_{K\in\pT_{h,\Gamma}}\Big(\gamma - \frac{\epsilon_c}{h_K}\Big) \| \llbracket u \rrbracket\|^2_{L^2(\Gamma_K)}                                                                                                                                                                       \\
           & \geqslant \sum_{i=1}^2 \|\alpha^{\half}\nabla u \|^2_{L^2(\Omega_i)} + \frac{1}{\epsilon_c}\| \curlii[]{\alpha \nabla_{\bn}u} \|_{H^{-\half}(\Gamma),h }
    - \sum_{K\in\pT_{h,\Gamma}} \Big(\frac{2C_\gamma}{\epsilon_c} \Big) \|\alpha^{\half} \nabla u\|_{L^2(K)}
    + \sum_{K\in\pT_{h,\Gamma}}\Big(\gamma - \frac{\epsilon_c}{h_K}\Big) \| \llbracket u  \rrbracket\|^2_{L^2(\Gamma_K)}                                                                                                                                                                      \\
           & = \sum_{K\in \widetilde{\pT}_{h}\setminus \pT_{h,\Gamma}} \|\alpha^{\half}\nabla u \|^2_{L^2(K)} + \sum_{K\in\pT_{h,\Gamma}} \|\alpha^{\half} \nabla u\|_{L^2(K)}
    + \frac{1}{\epsilon_c}\| \curlii[]{\alpha \nabla_{\bn}u} \|_{H^{-\half}(\Gamma),h } - \sum_{K\in\pT_{h,\Gamma}} \frac{2C_\gamma}{\epsilon_c} \|\alpha^{\half} \nabla u\|_{L^2(K)}                                                                                                         \\
           & \qquad+ \sum_{K\in\pT_{h,\Gamma}} \Big(\gamma - \frac{\epsilon_c}{h_K}\Big) \| \llbracket u \rrbracket\|^2_{L^2(\Gamma_K)}                                                                                                                                                       \\
           & = \sum_{K\in \widetilde{\pT}_{h}\setminus \pT_{h,\Gamma}} \|\alpha^{\half}\nabla u \|^2_{L^2(K)} + \sum_{K\in \pT_{h,\Gamma}} \half \|\alpha^{\half}\nabla u \|^2_{L^2(K)}
    + \frac{1}{\epsilon_c}\| \curlii[]{\alpha \nabla_{\bn}u} \|_{H^{-\half}(\Gamma),h }
    + \sum_{K\in\pT_{h,\Gamma}} \Big(\half - \frac{2C_\gamma}{\epsilon_c} \Big) \|\alpha^{\half} \nabla u\|_{L^2(K)}                                                                                                                                                                          \\
           & \qquad+ \sum_{K\in\pT_{h,\Gamma}} \Big(\gamma - \frac{\epsilon_c}{h_K}\Big) \| \llbracket u \rrbracket\|^2_{L^2(\Gamma_K)}.
  \end{align*}
}
{
  \begin{align*}
    a(u,u) & = \sum_{i=1}^2 \|\alpha^{\half}\nabla u \|^2_{L^2(\Omega_i)} - 2 (\curlii[]{\alpha \nabla_{\bn}u}, \llbracket u \rrbracket )_{L^2(\Gamma)}                                              \\
           & \quad  + \|\gamma^\half \llbracket u \rrbracket\|^2_{L^2(\Gamma)}                                                                                                                       \\
           & \geqslant  \sum_{i=1}^2 \|\alpha^{\half}\nabla u \|^2_{L^2(\Omega_i)} - \frac{1}{\epsilon_c} \| \curlii[]{\alpha \nabla_{\bn}u} \|_{H^{-\half}(\Gamma),h }                              \\
           & \quad -  \epsilon_c \|\llbracket u \rrbracket \|_{H^{\half}(\Gamma),h} + \gamma \|\llbracket u \rrbracket\|^2_{L^2(\Gamma)}                                                             \\
           & = \sum_{i=1}^2 \|\alpha^{\half}\nabla u \|^2_{L^2(\Omega_i)} + \Big(\frac{1}{\epsilon_c} - \frac{2}{\epsilon_c} \Big) \| \curlii[]{\alpha \nabla_{\bn}u} \|_{H^{-\half}(\Gamma),h }     \\
           & \quad + \sum_{K\in\pT_{h,\Gamma}}\Big(\gamma - \frac{\epsilon_c}{h_K}\Big) \| \llbracket u \rrbracket\|^2_{L^2(\Gamma_K)}                                                               \\
           & \geqslant \sum_{i=1}^2 \|\alpha^{\half}\nabla u \|^2_{L^2(\Omega_i)} + \frac{1}{\epsilon_c}\| \curlii[]{\alpha \nabla_{\bn}u} \|_{H^{-\half}(\Gamma),h }                                \\
           & \quad - \sum_{K\in\pT_{h,\Gamma}} \Big(\frac{2C_\gamma}{\epsilon_c} \Big) \|\alpha^{\half} \nabla u\|_{L^2(K)}                                                                          \\
           & \quad + \sum_{K\in\pT_{h,\Gamma}}\Big(\gamma - \frac{\epsilon_c}{h_K}\Big) \| \llbracket u  \rrbracket\|^2_{L^2(\Gamma_K)}                                                              \\
           & = \sum_{K\in \widetilde{\pT}_{h}\setminus \pT_{h,\Gamma}} \|\alpha^{\half}\nabla u \|^2_{L^2(K)} + \sum_{K\in\pT_{h,\Gamma}} \|\alpha^{\half} \nabla u\|_{L^2(K)}                       \\
           & \quad + \frac{1}{\epsilon_c}\| \curlii[]{\alpha \nabla_{\bn}u} \|_{H^{-\half}(\Gamma),h } - \sum_{K\in\pT_{h,\Gamma}} \frac{2C_\gamma}{\epsilon_c} \|\alpha^{\half} \nabla u\|_{L^2(K)} \\
           & \quad + \sum_{K\in\pT_{h,\Gamma}} \Big(\gamma - \frac{\epsilon_c}{h_K}\Big) \| \llbracket u \rrbracket\|^2_{L^2(\Gamma_K)}                                                              \\
           & = \sum_{K\in \widetilde{\pT}_{h}\setminus \pT_{h,\Gamma}} \|\alpha^{\half}\nabla u \|^2_{L^2(K)} + \sum_{K\in \pT_{h,\Gamma}} \half \|\alpha^{\half}\nabla u \|^2_{L^2(K)}              \\
           & \quad + \frac{1}{\epsilon_c}\| \curlii[]{\alpha \nabla_{\bn}u} \|_{H^{-\half}(\Gamma),h }                                                                                               \\
           & \quad + \sum_{K\in\pT_{h,\Gamma}} \Big(\half - \frac{2C_\gamma}{\epsilon_c} \Big) \|\alpha^{\half} \nabla u\|_{L^2(K)}                                                                  \\
           & \quad + \sum_{K\in\pT_{h,\Gamma}} \Big(\gamma - \frac{\epsilon_c}{h_K}\Big) \| \llbracket u \rrbracket\|^2_{L^2(\Gamma_K)}.
  \end{align*}
}
\endgroup
The second line uses Young's inequality for some $\epsilon_c > 0$, the fourth line follows the trace inequality \eqref{eq:inverse_ineq}.
 \end{appendix}

\end{document}